\documentclass[final,hidelinks,onefignum,onetabnum]{siamart220329}

\usepackage{amsmath,amsfonts}
\usepackage{algorithmic}
\usepackage{algorithm}
\usepackage{array}
\usepackage{textcomp}
\usepackage{stfloats}
\usepackage{url}
\usepackage[T1]{fontenc}
\usepackage{subcaption}
\usepackage{graphicx}
\usepackage{mathrsfs}
\usepackage[numbers,sort&compress]{natbib}
\usepackage{booktabs} 
\usepackage{multirow} 
\usepackage{amssymb}
\usepackage{xcolor}
\usepackage{hyperref}
\makeatletter
\IfFormatAtLeastTF{2024-06-01}{%
  \def\refstepcounter#1{\stepcounter{#1}%
    \def\@currentcounter{#1}%
    \protected@edef\@currentlabel
       {\csname p@#1\endcsname\csname the#1\endcsname}}%
}{}
\makeatother
\usepackage{cleveref}
\makeatletter
\def\cref@kernelcompat@noarg#1{%
  \cref@old@refstepcounter{#1}%
  \cref@constructprefix{#1}{\cref@result}%
  \@ifundefined{cref@#1@alias}%
    {\def\@tempa{#1}}%
    {\def\@tempa{\csname cref@#1@alias\endcsname}}%
  \protected@edef\cref@currentlabel{%
    [\@tempa][\arabic{#1}][\cref@result]%
    \csname p@#1\endcsname\csname the#1\endcsname}}
\def\cref@kernelcompat@optarg[#1]#2{%
  \cref@old@refstepcounter{#2}%
  \cref@constructprefix{#2}{\cref@result}%
  \@ifundefined{cref@#1@alias}%
    {\def\@tempa{#1}}%
    {\def\@tempa{\csname cref@#1@alias\endcsname}}%
  \protected@edef\cref@currentlabel{%
    [\@tempa][\arabic{#2}][\cref@result]%
    \csname p@#2\endcsname\csname the#2\endcsname}}
\AddToHook{begindocument}[zzz-kernelcompat]{%
  \IfFormatAtLeastTF{2024-06-01}{%
    \def\refstepcounter{%
      \@ifnextchar[{\refstepcounter@optarg}{\refstepcounter@noarg}}%
    \let\refstepcounter@noarg\cref@kernelcompat@noarg
    \let\refstepcounter@optarg\cref@kernelcompat@optarg
  }{}%
}
\makeatother

\usepackage{epstopdf}
\ifpdf
  \DeclareGraphicsExtensions{.eps,.pdf,.png,.jpg}
\else
  \DeclareGraphicsExtensions{.eps}
\fi

\newsiamremark{remark}{Remark}
\newsiamremark{hypothesis}{Hypothesis}
\crefname{hypothesis}{Hypothesis}{Hypotheses}
\newsiamthm{claim}{Claim}
\newsiamthm{assumption}{Assumption}
\crefname{assumption}{Assumption}{Assumption}
\headers{
Multistep Methods for Floquet Multipliers and Subspaces
}{Yehao Zhang, Yuncheng Xu, Chenyi Tan, and Yangfeng Su}
\title{Multistep Methods for Floquet Multipliers and Subspaces \thanks{%
\funding{This work was supported by the National Key R\&D Program of China under Grant Nos. 2020YFA0711900 and 2020YFA0711902.}}}

\author{Yehao Zhang\thanks{School of Mathematical Sciences, Fudan University, Shanghai, China 
  (\email{yehaozhang23@m.fudan.edu.cn}).}
  \and Yuncheng Xu\thanks{School of Mathematical Sciences, Fudan University, Shanghai, China 
  (\email{xuyc23@m.fudan.edu.cn}).}
  \and Chenyi Tan\thanks{School of Mathematical Sciences, Fudan University, Shanghai, China 
  (\email{cytan22@m.fudan.edu.cn}).}
\and Yangfeng Su \thanks{School of Mathematical Sciences, Fudan University, Shanghai, China 
  (\email{yfsu@fudan.edu.cn}).}}

\usepackage{amsopn}

\ifpdf
\hypersetup{
  pdftitle={Multistep Methods for Floquet Multipliers and Subspaces},
  pdfauthor={Yehao Zhang, Yuncheng Xu, Chenyi Tan, and Yangfeng Su}
}
\fi

\begin{document}
\maketitle
\begin{abstract}
 Accurate and efficient computation of Floquet multipliers and subspaces is essential for analyzing limit cycles in dynamical systems and periodic steady states in radio frequency circuit simulation. This problem is typically addressed by solving a periodic linear eigenvalue problem, which is obtained by discretizing the linear time-periodic system using one-step collocation methods. Collocation methods become costly for large-scale problems. Our alternative approach is to use multistep methods. A multistep method leads to a periodic polynomial eigenvalue problem (pPEP) and introduces additional parasitic periodic eigenvalues. We prove that, as the stepsize decreases, the computed Floquet multipliers and their associated invariant subspace converge at the consistency order, while the parasitic periodic eigenvalues converge to zero geometrically and hence become separated from the nonzero Floquet multipliers. We design a memory-efficient algorithm, pTOAR, to solve the large-scale pPEP. Its arithmetic and memory costs are almost independent of the choice of multistep methods when the pPEP arises from an implicit multistep discretization. Numerical results agree with our convergence analysis and demonstrate the efficiency of pTOAR.

\end{abstract}





\section{Introduction}
\label{sec:intro}
A linear time-periodic (LTP) system is described by
\begin{equation}
  \frac{\mathrm{d} x(t)}{\mathrm{d} t} = G(t) x(t),
  \label{equ:LTPOriginal}
\end{equation}
where \(G(t)\in \mathbb{R}^{n \times n}\) is a \(T\)-periodic matrix function, i.e., \(G(t + T) = G(t)\), and \(x(t)\in \mathbb{R}^{n}\) is the state vector. Let $X(t)$ denote the solution of the matrix initial-value problem $\dot X(t)=G(t)X(t)$ with $X(0)=I_n$. The monodromy matrix of \cref{equ:LTPOriginal} is $X(T)$; its eigenvalues and eigenvectors are the Floquet multipliers and their associated eigenvectors. 

Floquet multipliers and eigenvectors arise in dynamic simulation \cite{ALLEN2011DynamicalSimulation}, communications \cite{Gelli2003Communication}, control theory \cite{CHEN1997control,Bittanti1986LPSystem}, and radio frequency circuit simulation \cite{demir2000floquet}. In these applications, the stability of a limit cycle of a nonlinear autonomous system is determined by the Floquet multipliers of the LTP variational equation \cref{equ:LTPOriginal} evaluated along that cycle \cite{Lust2001improvedFloMul}.

The computation of Floquet multipliers and eigenvectors generally involves two key steps: discretization of \cref{equ:LTPOriginal} and solving the resulting eigenvalue problem. The simplest approach applies the backward Euler (BE) method, which yields a periodic eigenvalue problem (pEVP) solvable by the periodic QR algorithm \cite{bojanczyk1992Periodic, kressner2006PeriodicKS}. Existing high-precision packages such as \texttt{AUTO-07p} \cite{doedel2007auto} and \texttt{MATCONT} \cite{dhooge2008MATCONT} instead adopt high-order one-step collocation methods, which are well-suited for modest-scale problems with analytically available $G(t)$. For large-scale systems, however, collocation methods become expensive: they require extra function evaluations at interior collocation points, and the subsequent condensation step is costly and may destroy the inherent sparsity \cite{Lust2001improvedFloMul}.

In this paper we discretize \cref{equ:LTPOriginal} with multistep methods, which can achieve higher order than BE without requiring additional function evaluations. The price is that a $d$-step discretization induces a periodic polynomial eigenvalue problem (pPEP). Under the cyclic symmetry of the one-period discretization, this pPEP has $nd$ periodic eigenpairs and is therefore harder to solve than the one-step pEVP.

For a partition of $[0,T]$ into $p$ steps with maximum stepsize $h$, we prove that the pPEP spectrum separates into $n$ principal periodic eigenvalues converging to the Floquet multipliers and $(d-1)n$ parasitic periodic eigenvalues decaying geometrically to zero in the number of time steps $p$. For any spectrally separated cluster of $m$ Floquet multipliers, $m<n$, the associated periodic invariant subspace converges to its Floquet counterpart at order $\mathcal{O}(h^s)$, where $s$ is the consistency order of the adopted multistep method. The analysis rests on a uniformly conditioned moving basis in which the one-period recurrence is triangularized exactly by a periodic Riccati equation. This spectral separation permits an eigensolver to target only the principal periodic spectrum.

We propose \texttt{pTOAR}, an efficient solver for the pPEPs. It preserves the leading-order arithmetic complexity of periodic Krylov--Schur \cite{kressner2006PeriodicKS} applied to the companion linearization while avoiding the $d$-fold growth of the Arnoldi basis storage.

The remainder of this paper is organized as follows. In \cref{sec:pre}, the necessary background is reviewed. In \cref{sec:discrete}, we formulate the pPEP. In \cref{sec:Convergence}, we establish the convergence results. In \cref{sec:Algorithm}, we present the \texttt{pTOAR} algorithm and analyze its complexity. Numerical experiments are in \cref{sec:results}.

\section{Preliminaries}
\label{sec:pre}

\subsection{Floquet Theory}
\label{subsec:FloquetTheory}
Floquet theory provides the solution structure of the LTP system \cref{equ:LTPOriginal}; we recall the key results and refer to \cite{Amann1990ODEFLoquet} and \cite[Chapter~2]{farkas1994periodicmotions} for comprehensive treatments.

For the LTP system \cref{equ:LTPOriginal} with smooth, $T$-periodic $G(t)$, let $X(t)$ denote a \emph{fundamental solution matrix}, satisfying $\dot{X}(t) = G(t)X(t)$ and $X(0) = I_{n}$. The \emph{state transition matrix} $\Phi(t,s)$ maps the solution at time $s$ to that at time $t$, i.e., $\Phi(t,s) = X(t)X^{-1}(s).$
\begin{theorem}[Floquet Representation Theorem {\cite[Theorem 20.9]{Amann1990ODEFLoquet}}]
  \label{thm:FloquetRepresentation}
  For the LTP system \cref{equ:LTPOriginal}, there exist a constant matrix $M$ and a nonsingular $T$-periodic matrix function $U(t)\in C^{1}(\mathbb{R},\mathbb{C}^{n\times n})$ such that its state transition matrix $\Phi(t,s)$ can be written as
  \begin{equation}
    \Phi(t,s) = U(t) e^{M (t-s)} U(s)^{-1}.
    \label{equ:FloquetRepresentation}
  \end{equation}
\end{theorem}

Using the $T$-periodicity of $U(t)$, the one-period state transition matrix $\Phi(t+T,t)$ can be decomposed as
\begin{equation}
  \Phi(t+T,t) = U(t) \Lambda U(t)^{-1}.
  \label{equ:FloquetDecomposition}
\end{equation}
By \cref{equ:FloquetDecomposition}, all $\Phi(t+T,t)$ are similar to $\Lambda:=e^{MT}$ and share the same eigenvalues, which are called Floquet multipliers. Assume from now on that the multipliers are semisimple. Then, without loss of generality, we may choose \(\Lambda\) to be diagonal. Consequently, \cref{equ:FloquetDecomposition} is an eigendecomposition of \(\Phi(t+T,t)\). The columns of \(U(t)\) and the rows of \(U(t)^{-1}\) are the associated right and left Floquet eigenfunctions, respectively.

At $t=0$, \cref{equ:FloquetDecomposition} yields the Floquet eigenvalue problem \cref{equ:FloquetEVP}:
\begin{equation} 
    \Phi(T,0)u(0) = u(0)\lambda.
  \label{equ:FloquetEVP}
\end{equation}
Here, $\lambda$ is a Floquet multiplier, and $u(0)$ is the \emph{Floquet eigenvector} associated with the multiplier; the corresponding periodic eigenfunction $u(t)$ can be recovered from \cref{equ:FloquetRepresentation}.

\subsection{One-step Methods for Solving Floquet Eigenpairs}
\label{subsec:OneStep}
Let the period $[0,T]$ be partitioned at time points $0 = t_0 < t_1 < \cdots < t_p = T$, with stepsizes $h_{i} = t_{i}-t_{i-1}$ allowed to be variable, and with maximum $h = \max_{i = 1}^{p}h_{i}$. Rewrite \cref{equ:FloquetEVP} on the time grid:
\[\Phi(t_{p},t_{p-1})\cdots\Phi(t_{2},t_{1})\Phi(t_{1},t_{0})u(0) = u(0)\lambda.\]
When $h$ is small enough, one-step methods approximate each $\Phi(t_{i},t_{i-1})$ by a matrix $A^{(i)}$ for $i=1,\ldots,p$, and approximate the eigenpair $\left(u(0),\lambda\right)$ by  $\left(u_{0},\hat{\lambda}\right)$. This leads to a periodic eigenvalue problem \cite{kressner2006PeriodicKS}:
\begin{equation}
  A^{(p)}\cdots A^{(2)}A^{(1)} u_{0} = u_{0}\hat{\lambda}.
  \label{equ:OneStepEVP}
\end{equation}

One of the simplest one-step methods is the backward Euler (BE) method. BE requires only one function evaluation of $G(t)$ at time node $t_{i}$ on each interval $[t_{i-1},t_{i}]$; when BE is applied to \cref{equ:LTPOriginal}, $A^{(i)} = \left[I-G(t_{i})h_{i}\right]^{-1}$. However, its first-order accuracy is usually insufficient in practice.

High-order one-step methods, such as collocation, improve discretization accuracy. In what follows, we briefly outline the collocation discretization; a comprehensive treatment can be found in \cite{Lust2001improvedFloMul}.

On each interval $[t_{i-1},t_i]$, collocation methods introduce $d$ internal points to construct a polynomial approximation to $x(t)$. Then, the interpolant is required to satisfy \cref{equ:LTPOriginal} at $d$ Gaussian collocation points within $[t_{i-1},t_{i}]$. Consequently, the system matrix $G(t)$ must be evaluated at these internal points in addition to the prescribed time-grid points.

The next step, known as \emph{condensation}, eliminates the internal point variables to obtain the one-step relation \(x_{i} = A^{(i)} x_{i-1}\), typically via an LU elimination \cite{Lust2001improvedFloMul}. This process yields an explicitly formed, typically dense local state transition matrix $A^{(i)}$.

Collocation methods are effective when $G(t)$ is analytically available and the problem is moderate-size, and they are used in established packages such as \texttt{MATCONT} \cite{dhooge2008MATCONT} and \texttt{AUTO-07p} \cite{doedel2007auto}. For large-scale systems, or when $G(t)$ is available only on a prescribed time grid, however, the required internal evaluations, condensation cost, and loss of sparsity make collocation methods expensive or even inapplicable.

\subsection{Multistep Methods for Solving Initial Value Problems} 
Multistep methods offer an alternative approach to enhancing discretization accuracy by using previously computed values at neighboring grid points. We recall the variable stepsize form, consistency, and zero-stability needed below \cite{hairer1993Multistep}.

\begin{definition}[Variable Stepsize Multistep Methods {\cite[Chapter 3.5]{hairer1993Multistep}}]
  Let $0= t_{0}<\cdots<t_{p} = T$ be a time grid on the interval $[0,T]$. We define the stepsize $h_{i} = t_{i} - t_{i-1}$ and the stepsize ratio $\omega_{i} = h_{i}/h_{i-1}$. The ratio vector $\boldsymbol{\omega}_i \in \mathbb{R}^{d-1}$ is defined by \(\boldsymbol{\omega}_i = (\omega_{i-d+2},\ldots,\omega_{i}).\)  For an ODE $\dot{x} = f(t,x)$, a $d$-step method updates the solution at $t_{i}$ according to:
  \begin{equation}
    \sum_{j=0}^{d}\alpha_{j,i} x_{i-d+j} = h_{i} \sum_{j=0}^{d}\beta_{j,i} f(t_{i-d+j},x_{i-d+j}),\quad i = d,\ldots,p,
    \label{def:MultistepMethods}
  \end{equation}
  where the coefficients are determined by smooth functions of the local step ratios, i.e., $\alpha_{j,i} = \alpha_{j}(\boldsymbol{\omega}_i)$ and $\beta_{j,i} = \beta_{j}(\boldsymbol{\omega}_i)$, with $\alpha_{d,i}\neq 0$ and $\left\lvert \alpha_{0,i} \right\rvert +\left\lvert \beta_{0,i} \right\rvert \neq 0$, which is the irreducibility condition.
\end{definition}

For detailed derivations of these coefficients, see \cite[Chapter 3]{hairer1993Multistep}. As an illustration, the second-order backward differentiation formula (BDF2, also called Gear2) has coefficients given by \cref{def:VariableStepSizeGear2}. 
\begin{equation}
\alpha_{0,i} = \frac{{\omega_{i}}^{2}}{1+2\omega_{i}},\,\alpha_{1,i} = -1-\alpha_{0,i},\alpha_{2,i} = 1,\,\beta_{2,i} = \frac{1+\omega_{i}}{1+2\omega_{i}},\, \beta_{1,i} = \beta_{0,i} = 0.
\label{def:VariableStepSizeGear2}
\end{equation}

We next recall consistency and stability for such multistep methods.
\begin{definition}[Consistency {\cite[Chapter~III.5, Eq.(5.17)]{hairer1993Multistep}}]
  \label{def:Consistency}
  A multistep method \cref{def:MultistepMethods} is consistent of order $s$ if, when applied to a sufficiently smooth function $x(t)$, its local truncation error $R_{i}$ satisfies $\left\lVert R_{i} \right\rVert  = \mathcal{O}(h_{i}^{s+1})$ as $h_{i} \to 0$, where
  \begin{equation*}
  R_{i} = \sum_{j=0}^{d}\alpha_{j,i}x(t_{i-d+j})-h_{i}\sum_{j=0}^{d}\beta_{j,i}\frac{\mathrm{d}x}{\mathrm{d}t}(t_{i-d+j}).
  \end{equation*}

\end{definition}
 
The stability of a multistep method guarantees the boundedness of the numerical solution. We first recall zero-stability (also called Dahlquist stability), which is defined for constant stepsizes.
\begin{definition}[Zero-stability \cite{dahlquist1956Convergence}]
  \label{def:ZeroStability}
  A $d$-step method is called zero-stable if its first generating polynomial 
  \[\rho(\nu)\equiv \sum_{j=0}^{d} \alpha_{j}(\mathbf{1}) \nu^{j},\mathbf{1} = (1,\ldots,1)\in\mathbb R^{d-1}\] satisfies the root condition:
  \begin{enumerate}
    \item All roots of $\rho(\nu)$ lie on or within the unit circle;
    \item The roots on the unit circle are simple.
  \end{enumerate}
\end{definition}
If the only root of $\rho(\nu)$ on the unit circle is $\nu =1$, then the method is called \emph{strongly zero-stable}.
For a consistent multistep method, $\nu_1=1$ is called the
\emph{principal root}; all other roots $\nu_2,\ldots,\nu_d$, counted with
algebraic multiplicity, are called the \emph{parasitic roots}. Thus, for a
strongly zero-stable method, $|\nu_\tau|<1$, $\tau=2,\ldots,d$.

\section{Multistep Discretization of an LTP System and Derived pPEP}
\label{sec:discrete}
Applying a multistep method to \cref{equ:LTPOriginal} yields a \emph{periodic difference equation}, which can be formulated as a \emph{periodic polynomial eigenvalue problem} (pPEP); linearizing the pPEP then yields an equivalent \emph{periodic eigenvalue problem} (pEVP).

\subsection{Multistep Discretization of an LTP System Induces a pPEP}
Consider a partition of one period $[0,T]$ at time points $0 = t_0 < t_1 < \cdots < t_p = T$. This grid is extended to an infinite sequence $\left\{t_{i}\right\}_{i \in \mathbb{Z}}$ satisfying $t_{i+p} = t_{i} + T$, which defines a \emph{periodic time grid}. A quantity that is $p$-periodic on this grid carries the phase superscript $(i)$. Accordingly, the multistep coefficients $\left(\alpha_{j,i},\beta_{j,i}\right)$ are denoted by $\left(\alpha_{j}^{(i)},\beta_{j}^{(i)}\right)$, and the stepsizes and stepsize ratios are written as $h^{(i)} = t_{i}-t_{i-1}$ and $\omega^{(i)}=h^{(i)}/h^{(i-1)}$, respectively. We also set the maximum step $h= \max_{i=1}^{p}h^{(i)}$, and the ratio vector $\boldsymbol{\omega}^{(i)} := (\omega^{(i-d+2)},\ldots,\omega^{(i)}) \in \mathbb{R}^{d-1}$.

Applying a linear $d$-step method to \cref{equ:LTPOriginal} yields a linear difference equation of order $d$ of the form
\begin{equation}
\sum_{j=0}^{d} A_j^{(i)} x_{i-d+j} = 0,
\label{equ:kstepRecursion}
\end{equation}
where the coefficient matrices $A_j^{(i)} \in \mathbb{R}^{n\times n}$ are $p$-periodic in $i$ and given by
\begin{equation}
A_j^{(i)} = \alpha_j^{(i)} I_n - \beta_j^{(i)} G(t_{i-d+j}) h^{(i)}, 
\quad j=0,\ldots,d.
  \label{def:CoefficientMat}
\end{equation}

To compute Floquet multipliers, we seek an eigenvalue formulation associated with \cref{equ:kstepRecursion}. Motivated by the discrete Floquet theorem for one-step periodic recursions \cite{DACUNHA2011UnifiedFloquet}, we substitute the Floquet-type ansatz $x_i = u^{(i)} \theta^{i}$ into the general $d$-step recursion \cref{equ:kstepRecursion}, where $u^{(i)}$ is $p$-periodic.

\begin{definition}[Periodic Polynomial Eigenvalue Problem (pPEP)]
\label{def:periodicPEP}\\
Let $\{A_j^{(i)}\}_{j=0,\ldots,d}^{i=1,\ldots,p}$ denote $d+1$ sequences of matrices in $\mathbb{R}^{n\times n}$, each of period $p$. The \emph{periodic polynomial eigenvalue problem (pPEP)} consists of finding a scalar $\theta \in \mathbb{C}$ and a nonzero $p$-periodic sequence $\{u^{(i)}\}_{i\in\mathbb{Z}}$,  $u^{(i)}\in \mathbb{C}^n$ such that
\begin{equation}
  \sum_{j=0}^{d} A_{j}^{(i)}\, u^{(i-d+j)} \theta^{j} = 0,
  \qquad i = 1,\ldots,p,
  \label{equ:pPEPDefinition}
\end{equation}
where all parenthesized phase superscripts $(i)$ are understood cyclically with period $p$.
\end{definition}
Stacking $p$ components in \cref{equ:pPEPDefinition} leads to a degree-$d$ polynomial eigenvalue problem (PEP) \cref{equ:EVPPolynomialMatrix} \cite{mackey2015PEP} with block-cyclic coefficient matrices. To illustrate this, we define the block-diagonal matrices \(\mathscr{A}_j = \operatorname{diag}_{i=1}^{p}(A_j^{(i)}),\, j = 0, \ldots, d,\) and the cyclic shift operator \(\mathscr{C} = C_p \otimes I_n,\,C_p = [e_2, \dots, e_p, e_1],\) where $e_i$ denotes the $i$-th standard basis vector of $\mathbb{R}^p$.
\begin{equation}
  \left( \sum_{j=0}^{d} \mathscr{A}_j \mathscr{C}^{d-j} \theta^j \right) \operatorname{col}(u^{(i)})_{i=1}^{p} = 0,\quad \operatorname{col}(u^{(i)})_{i=1}^{p} \equiv \begin{bmatrix} u^{(1)} \\ u^{(2)} \\ \vdots \\ u^{(p)} \end{bmatrix} \in \mathbb{C}^{np}.
  \label{equ:EVPPolynomialMatrix}
\end{equation}

If $(\theta, \operatorname{col}\{u^{(i)}\}_{i=1}^{p})$ is an eigenpair, then for each $k=0,\ldots,p-1$,
\[
\left(\theta \zeta_p^{k}, \operatorname{col}\{u^{(i)} \zeta_p^{-ik}\}_{i=1}^{p}\right)
\]
is also an eigenpair, where $\zeta_p = e^{\jmath 2\pi / p}$ is a primitive $p$-th root of unity with $\jmath^2=-1$. For a nonzero root orbit, choose the representative $\theta_{\star}$ with $\arg(\theta_{\star})\in[0,2\pi/p)$. If $\theta=0$, take $\theta_{\star} = 0$ and retain its associated periodic sequence. In either case, define the \emph{periodic eigenvalue} by $\lambda=\theta_{\star}^p$ and the \emph{periodic eigenvector} by the sequence associated with $\theta_{\star}$.

The pPEP is called \emph{monic} if its leading coefficient matrices satisfy $A_d^{(i)} = I_n$ for all $i=1,\ldots,p$. For sufficiently small maximum stepsize $h$, the matrices $A_d^{(i)}$ defined in \cref{def:CoefficientMat} are nonsingular with uniformly bounded inverse, and left-multiplying the $i$-th equation in \cref{equ:pPEPDefinition} by $\left(A_d^{(i)}\right)^{-1}$ yields an equivalent monic pPEP. Thus, we denote
\[\widetilde A_j^{(i)}:=(A_d^{(i)})^{-1}A_j^{(i)},\qquad
\widetilde{\mathscr{A}}_j:=\operatorname{diag}_{i=1}^{p}
(\widetilde A_j^{(i)}),\quad j=0,\ldots,d-1.\]

\subsection{Linearized pPEP as a Periodic Eigenvalue Problem}
After monic normalization, we transform \cref{equ:EVPPolynomialMatrix} into a pEVP of dimension $nd$ and period $p$ in two steps: a companion linearization followed by a permutation.
\begin{equation}
  \begin{bmatrix}
    &I_{np}&&\\
    &&\ddots&\\
    &&&I_{np}\\
    -\widetilde{\mathscr{A}}_{0}\mathscr{C}^{d}&-\widetilde{\mathscr{A}}_{1}\mathscr{C}^{d-1}&\cdots&-\widetilde{\mathscr{A}}_{d-1}\mathscr{C}\\
  \end{bmatrix}  
  \begin{bmatrix}
    \operatorname{col}(u^{(i)})_{i=1}^{p}\\
    \operatorname{col}(u^{(i)})_{i=1}^{p}\theta\\
    \vdots\\
    \operatorname{col}(u^{(i)})_{i=1}^{p}\theta^{d-1}
  \end{bmatrix}=  \begin{bmatrix}
    \operatorname{col}(u^{(i)})_{i=1}^{p}\\
    \operatorname{col}(u^{(i)})_{i=1}^{p}\theta\\
    \vdots\\
    \operatorname{col}(u^{(i)})_{i=1}^{p}\theta^{d-1}
  \end{bmatrix}\theta.
  \label{equ:EVPLinearizedPEP}
\end{equation}

The eigenvector in \cref{equ:EVPLinearizedPEP} exhibits a two-level structure: the inner structure encodes the periodicity of the coefficients, while the outer structure arises from the linearization of \cref{equ:EVPPolynomialMatrix}. To make this structure explicit, we introduce a permutation matrix $\mathscr{P}$ that reorders the stacked vector:
  \[\mathscr{P}\begin{bmatrix}
    \operatorname{col}(u^{(i)})_{i=1}^{p}\\
    \operatorname{col}(u^{(i)})_{i=1}^{p}\theta\\
    \vdots\\
    \operatorname{col}(u^{(i)})_{i=1}^{p}\theta^{d-1}
  \end{bmatrix}=\begin{bmatrix}
  \operatorname{col}(u^{(1-d+j)}\theta^{j-1})_{j=1}^{d}\\
  \operatorname{col}(u^{(2-d+j)}\theta^{j-1})_{j=1}^{d}\\
  \vdots\\
  \operatorname{col}(u^{(p-d+j)}\theta^{j-1})_{j=1}^{d}
  \end{bmatrix},\, \operatorname{col}(u^{(i-d+j)}\theta^{j-1})_{j=1}^{d} = \begin{bmatrix}
  u^{(i-d+1)}\\u^{(i-d+2)}\theta\\\vdots\\u^{(i)}\theta^{d-1}\\
\end{bmatrix}.
  \]
After this permutation, \cref{equ:EVPLinearizedPEP} takes the following block-cyclic form:
\begin{equation}
  \begin{bmatrix}
    &&&\mathcal{L}^{(1)}\\
    \mathcal{L}^{(2)}&&&\\
    &\ddots&&\\
    &&\mathcal{L}^{(p)}&\\
  \end{bmatrix}\begin{bmatrix}
  \operatorname{col}(u^{(1-d+j)}\theta^{j-1})_{j=1}^{d}\\\operatorname{col}(u^{(2-d+j)}\theta^{j-1})_{j=1}^{d} \\\vdots\\\operatorname{col}(u^{(p-d+j)}\theta^{j-1})_{j=1}^{d}\\
\end{bmatrix}=  \begin{bmatrix}
  \operatorname{col}(u^{(1-d+j)}\theta^{j-1})_{j=1}^{d}\\\operatorname{col}(u^{(2-d+j)}\theta^{j-1})_{j=1}^{d} \\\vdots\\\operatorname{col}(u^{(p-d+j)}\theta^{j-1})_{j=1}^{d}\\
\end{bmatrix}\theta,
  \label{equ:EVPLinearizedCyclic} 
\end{equation}
where each block $\mathcal{L}^{(i)}$ is the companion matrix of \cref{equ:kstepRecursion} at step $i$. 
\begin{equation}
  \mathcal{L}^{(i)} = \begin{bmatrix}
    &I_{n}&&\\
    &&\ddots&\\
    &&&I_{n}\\
    -\widetilde A_{0}^{(i)}&-\widetilde A_{1}^{(i)}&\cdots&-\widetilde A_{d-1}^{(i)}\\
  \end{bmatrix}.
    \label{def:L}
\end{equation}

The block-cyclic structure implies that \cref{equ:EVPLinearizedCyclic} reduces to $p$ periodic eigenvalue problems of dimension $nd$, one for each starting index $i = 1,\ldots,p$:
\begin{equation}
\mathcal{L}^{(i+p)}\mathcal{L}^{(i+p-1)}\cdots \mathcal{L}^{(i+1)}  \operatorname{col}(u^{(i-d+j)}\theta^{j-1})_{j=1}^{d} =  \operatorname{col}(u^{(i-d+j)}\theta^{j-1})_{j=1}^{d}\theta^{p},\,i = 1,\ldots,p.
\label{equ:EVPLinearizedPeriodic}
\end{equation}
Since linearization and permutation are invertible, \cref{equ:EVPLinearizedCyclic} and \cref{equ:EVPLinearizedPeriodic} are equivalent to the original pPEP \cref{equ:pPEPDefinition} if all $A_{d}^{(i)}$ are nonsingular.

\section{Convergence Analysis}
\label{sec:Convergence}
This section establishes how the linearized pPEP \cref{equ:EVPLinearizedPeriodic} relates to the continuous Floquet problem \cref{equ:FloquetEVP}. Each companion matrix $\mathcal{L}^{(i)}$ in \cref{def:L} advances the $d$-stacked state by one step \cite{gohberg2005matrixpolynomials}; hence the propagator from step $i$ to $i+j$ is $\mathcal{F}_{i+j,i}=\mathcal{L}^{(i+j)}\cdots \mathcal{L}^{(i+1)}$. Since $\mathcal{L}^{(i)}$ is $p$-periodic and cyclic shifts of square factors preserve the spectrum, all one-period propagators $\mathcal{F}_{i+p,i}$ share the same spectrum; it therefore suffices to consider $\mathcal{F}_{p,0}$. By \cref{equ:EVPLinearizedPeriodic}, the periodic eigenvalues of the pPEP are exactly the eigenvalues of $\mathcal{F}_{p,0}$.

The Floquet eigenvalue problem \cref{equ:FloquetEVP} characterizes solutions of \cref{equ:LTPOriginal} with $x(T)=\lambda x(0)$ and $x(0)=u(0)$. Thus, its $n$ multipliers provide the principal reference spectrum of $\mathcal{F}_{p,0}$; the remaining $(d-1)n$ periodic eigenvalues arise from the augmented multistep state and are called parasitic periodic eigenvalues \cite{citro2020parasitism,dambrosio2014parasitism,hairer1999backwarderrorMultistep}.

\subsection{Standing Assumptions and the Main Convergence Result}
\label{subsec:AssumptionsMainResult}

We consider a fixed $d$-step method with $d\ge2$. Unless the subscript $F$ is shown, $\left\lVert \cdot \right\rVert$ denotes the spectral $2$-norm; $\left\lVert \cdot \right\rVert_{\mathrm{F}}$ denotes the Frobenius norm. All constants are independent of the index $i$, of $p$, and of the admissible periodic grids as the maximum stepsize $h\to0$. We impose the following assumptions.

\begin{enumerate}
  \item[\textbf{(A1)}] \textbf{Regularity and uniform smoothness.}
  The periodic matrix $G$ in \cref{equ:LTPOriginal} satisfies $G \in C^{s}_{\mathrm{per}}([0,T],\mathbb{R}^{n\times n})$, so every solution of \cref{equ:LTPOriginal} is $C^{s+1}$; in particular $\left\lVert G(t) \right\rVert\leq C_G$.

  \item[\textbf{(A2)}] \textbf{Admissible periodic grids and smooth coefficients.}
  The grids are periodic, $t_{i+p} = t_i + T$, satisfy
  \[
    \sum_{i=1}^p h^{(i)}=T,\qquad
    h=\max_{1\leq i\leq p} h^{(i)}\to0,\qquad
    \sum_{i=1}^p|\omega^{(i)}-1|\leq C_T,
  \]
  and their ratio vectors lie in a common compact set
  $K_\omega\subset\mathbb R^{d-1}$ containing $\mathbf{1}$ and every line segment joining $\mathbf{1}$ to an admissible ratio vector. The coefficient functions satisfy $\alpha_j,\beta_j\in C^1(K_\omega)$ with uniform bounds, and $\alpha_d$ is bounded away from zero on $K_\omega$.

  \item[\textbf{(A3)}] \textbf{Uniform consistency.}
  The method has consistency order $s\ge1$ in the sense of \cref{def:Consistency}, uniformly over the admissible grid.
  \item[\textbf{(A4)}] \textbf{Strong zero-stability with simple nonzero parasitic roots.} The method is strongly zero-stable in the sense of \cref{def:ZeroStability}, and its parasitic roots are distinct and nonzero.
  \item[\textbf{(A5)}] \textbf{Semisimple Floquet spectrum.}
  The monodromy matrix $\Phi(T,0)$ of \cref{equ:LTPOriginal} is diagonalizable over $\mathbb C$; equivalently, all Floquet multipliers are semisimple.
\end{enumerate}

The assumptions have distinct roles. Assumptions~(A1), (A2), (A3) together with strong zero-stability in (A4) are standard assumptions for multistep discretization. The distinct and nonzero root condition in (A4) and semisimplicity in (A5) support the moving-basis and invariant-subspace arguments below. Repeated parasitic roots and nonsemisimple Floquet multipliers require separate analyses and are outside the present scope.

To state the convergence result, choose an eigenvector matrix $U(0)$ by \textup{(A5)}, such that
\begin{equation*}
  \Phi(T,0)U(0)=U(0)\Lambda,
  \qquad
  \Lambda=\operatorname{diag}\left(\lambda_1,\ldots,\lambda_n\right).
\end{equation*}

For $1\leq m<n$, order the selected multipliers $\Lambda_{m}:=\operatorname{diag}\left(\lambda_1,\ldots,\lambda_{m}\right)$ first and write
\[
  U(0)=\begin{bmatrix}U_m(0)&U_{n-m}(0)\end{bmatrix},
  \qquad
  \Lambda=\operatorname{diag}(\Lambda_m,\Lambda_{n-m}).
\]
The selected block may contain repeated multipliers; only its separation from the complementary block is required. We measure the separation between the selected and complementary Floquet subspaces by
\begin{equation}
  \gamma_{m}
  :=\inf_{0\ne Z\in\mathbb C^{(n-m)\times m}}
  \frac{\left\lVert \Lambda_{n-m}Z-Z\Lambda_m \right\rVert_{\mathrm{F}}}{\left\lVert Z \right\rVert_{\mathrm{F}}}
  =\min_{\substack{\lambda\in\sigma(\Lambda_m)\\
                    \zeta\in\sigma(\Lambda_{n-m})}}
    |\lambda-\zeta|.
  \label{def:ClusterSep}
\end{equation}
The equality follows because both blocks are diagonal \cite{Stewart1971ISbounds}. For two \(m\)-dimensional subspaces \(\mathcal{X},\mathcal{Y}\subset\mathbb{C}^{N}\), let \(Q_{\mathcal{X}},Q_{\mathcal{Y}}\in\mathbb{C}^{N\times m}\) be orthonormal basis matrices for \(\mathcal{X}\) and \(\mathcal{Y}\), respectively. The canonical angles \(0\leq\theta_{1}\leq\cdots\leq\theta_{m}\leq\pi/2\) are defined by \(\sigma_i\!\left(Q_{\mathcal{X}}^{\mathrm H}Q_{\mathcal{Y}}\right) =\cos\theta_i,\, i=1,\ldots,m, \) where the singular values are ordered nonincreasingly. We collect the canonical angles in the diagonal matrix \(\Theta(\mathcal{X},\mathcal{Y})
:=\operatorname{diag}(\theta_{1},\ldots,\theta_{m}).\)

For $1\leq m< n$, define the exact stacked Floquet subspace
\begin{equation}
  \mathcal{U}_m:=\operatorname{range}\left(
    \operatorname{col}\left(
      \Phi(t_{-d+1},0)U_m(0),\ldots,
      \Phi(t_{0},0)U_m(0)
    \right)
  \right)
  \subset\mathbb C^{dn}.
  \label{def:ExactStackedSelectedModes}
\end{equation}
The dependence on the grid enters only through the points $t_{-d+1},\ldots,t_{0}$, which are determined by the admissible grid family as $h\to0$.

\begin{theorem}[Convergence of Floquet multipliers and invariant subspaces]
\label{thm:MainConvergence} 

Consider the pPEP obtained by applying a $d$-step multistep method of consistency order $s$ to \cref{equ:LTPOriginal}, together with its companion linearization. Let $\mathcal{F}_{p,0}$ denote the one-period propagator of the companion linearization. Suppose that assumptions \textup{(A1)}--\textup{(A5)} hold. Along any admissible periodic-grid family covered by \textup{(A2)}, as the maximum stepsize \(h\to0\), and hence, \(p\to\infty,\) the $nd$ periodic eigenvalues \(\sigma(\mathcal{F}_{p,0})=\{\widehat{\lambda}_{k}\}_{k=1}^{nd},\) counted with algebraic multiplicity, split into $n$ principal periodic eigenvalues and $(d-1)n$ parasitic periodic eigenvalues with the following properties.
\begin{enumerate}
\item After suitable ordering, the $n$ principal periodic eigenvalues satisfy
\[
  \widehat{\lambda}_{k}=\lambda_{k}+\mathcal{O}(h^s),
  \qquad k=1,\ldots,n.
\]
\item Let $r := \max_{2\le\tau\le d}|\nu_{\tau}|$. There exists a constant $C>0$, independent of the grid, such that the parasitic periodic eigenvalues satisfy
$$
|\widehat{\lambda}_{k}|\leq Cr^p,\qquad k=n+1,\ldots,nd.
$$

\item Given $m<n$ selected Floquet multipliers $\lambda_1,\ldots,\lambda_m$ whose separation from the complementary spectrum satisfies $\gamma_m>0$, the corresponding principal periodic eigenvalues define an invariant subspace $\widehat{\mathcal{U}}_m$ of $\mathcal{F}_{p,0}$,  satisfying
\[
\left\lVert \sin\Theta\left(\widehat{\mathcal{U}}_m,\mathcal{U}_m \right) \right\rVert 
  =\mathcal{O}\left(\frac{h^s}{\gamma_{m}}\right).
\]
\end{enumerate}
\end{theorem}

Part~(1) transfers the consistency order of the discretization to the principal periodic eigenvalues. Part~(2) shows that the parasitic spectrum decays exponentially with the number of time steps and therefore separates from the nonzero Floquet multipliers on sufficiently fine grids. When Part~(3) is specialized to \(m=1\), \cref{def:ExactStackedSelectedModes} shows that \(\mathcal{U}_1\) is spanned by the \(d\)-point starting history \(\operatorname{col}\left(\Phi(t_{-d+1},0)u(0),\ldots,\Phi(t_{0},0)u(0)\right)\) of the ODE solution associated with the target Floquet multiplier \(\lambda\). Part~(3) therefore provides, up to a scalar normalization, an \(\mathcal{O}(h^s/\gamma_1)\) approximation to this starting history. These approximate values can be used to initialize the multistep recurrence and compute the solution at the remaining grid points over one period.

The proof consists of three stages. First, we introduce a uniformly conditioned moving basis $\mathcal{W}_{i}$, in which each companion step satisfies $\mathcal{W}_{i}^{-1}\mathcal{L}^{(i)}\mathcal{W}_{i-1}=\mathcal{D}+\mathcal{E}^{(i)}$, where $\mathcal{D}$ is the reference diagonal transition and $\mathcal{E}^{(i)}$ will be shown to be small. Second, a periodic Riccati equation eliminates the lower-left block and exactly separates the principal and parasitic products. Finally, one-period estimates prove the principal and parasitic eigenvalue bounds, while an invariant-subspace perturbation argument within the principal block proves part~(3).

\subsection{Moving Basis and Local Transition Bounds}
\label{subsec:MovingBasis}
We combine the state-transition histories with the constant parasitic modes. The resulting basis is uniformly conditioned and satisfies the Floquet quasi-periodicity relation used in the Riccati boundary condition.

For the principal part, write $\Phi(t):=\Phi(t,0)$ and
\[
W_{\mathrm{pr},i}:=\operatorname{col}\left(\Phi(t_{i-d+1}),\ldots,\Phi(t_i)\right).
\]
For each parasitic root $\nu_{\tau}$ of the first generating polynomial $\sum_{j=0}^d\alpha_j(\mathbf1)\nu^j$ (see \cref{def:ZeroStability}),  introduce the associated constant block column
\[
v_{\tau} :=
\begin{bmatrix}
1 & \nu_{\tau} & \cdots & \nu_{\tau}^{d-1}
\end{bmatrix}^{\mathrm T}\otimes I_n
\in \mathbb{C}^{dn\times n},\quad \tau=2,\ldots,d.
\]
The columns $v_2,\ldots,v_d$ constitute the constant parasitic part of the moving basis.
\begin{definition}[Moving basis]
  \label{def:MovingBasisW}
The \emph{moving basis} at the grid point $t_i$ is
\begin{equation*}
  \begin{aligned}
      \mathcal{W}_{i} :=&\begin{bmatrix}
        \Phi(t_{i-d+1}) & I_{n} & \cdots & I_{n}\\
        \Phi(t_{i-d+2}) & \nu_{2}I_{n} &\cdots &\nu_{d}I_{n}\\
        \vdots & \vdots & \ddots & \vdots\\
        \Phi(t_{i}) & \nu_{2}^{d-1}I_{n} & \cdots & \nu_{d}^{d-1}I_{n}
      \end{bmatrix}
  \end{aligned}
\end{equation*}
\end{definition}

\begin{lemma}[Uniform conditioning of the moving basis]
  \label{lem:MovingBasisCondition}
  Under \textup{(A1)}, \textup{(A2)}, and \textup{(A4)}, there exist $h_{0}>0$ and $C_W>0$, independent of $i$ and $p$, such that for all $h\leq h_{0}$,
  \[
    \sup_{0\leq i\leq p} \left(\left\lVert \mathcal{W}_{i} \right\rVert + \left\lVert \mathcal{W}_{i}^{-1} \right\rVert\right) \leq C_W .
  \]
\end{lemma}
\begin{proof}
The key observation is that $\mathcal{W}_{i}$ is a perturbation of a fixed invertible block Vandermonde matrix. Since $\Phi(t_j)=\Phi(t_j,t_i)\Phi(t_i)$, we may factor
\[
  \mathcal{W}_{i} = \mathcal{V}(t_i)\,\operatorname{diag}\left(\Phi(t_i), I_{(d-1)n}\right),
\]
where the first block column of $\mathcal{V}(t_i)$ is $\operatorname{col}(\Phi(t_{i-d+1},t_i),\ldots,\Phi(t_{i-1},t_i),I_n)$ and the remaining block columns are $v_2,\ldots,v_d$. Since $|t_{i-k}-t_i|\leq (d-1)h$ and $\left\lVert \Phi(t_j,t_i) - I_n \right\rVert \leq C_V h$ by (A1), $\mathcal{V}(t_i)$ differs from the block Vandermonde matrix $\mathcal{V}_* = V\otimes I_n$, $V_{k\tau} = \nu_{\tau}^{k-1}$, by $\mathcal{O}(h)$ uniformly in $i$. The nodes $1,\nu_2,\ldots,\nu_d$ are distinct by Assumption (A4), so $V$, and hence $\mathcal{V}_*$, is nonsingular. For $h_{0}$ with $\left\lVert \mathcal{V}_*^{-1} \right\rVert\,\left\lVert \mathcal{V}(t_i)-\mathcal{V}_* \right\rVert \leq \tfrac12$, the Neumann series gives $\left\lVert \mathcal{V}(t_i)^{-1} \right\rVert \leq 2\left\lVert \mathcal{V}_*^{-1} \right\rVert$. Finally, $\Phi(t)$ and $\Phi(t)^{-1}$ are continuous on the compact interval $[0,T]$ by assumption (A1), hence uniformly bounded.
\end{proof}

By the periodicity of $G(t)$, the principal stacked block vector satisfies $W_{\mathrm{pr},i+p}=W_{\mathrm{pr},i}\Phi(T,0)$. Hence the moving basis satisfies the Floquet quasi-periodicity relation
\begin{equation}
  \mathcal{W}_{p}=\mathcal{W}_{0}\mathcal{B},\quad \mathcal{B} := \operatorname{diag}\left(\Phi(T,0), I_{(d-1)n}\right).
  \label{equ:MovingFrameQuasiPeriodicity}
\end{equation}
\begin{definition}[Local Transition Representation]
  \label{def:LocalTransitionPerturbation}
We write the local transition of $\mathcal{L}^{(i)}$ in the moving basis $\left\{\mathcal{W}_{i}\right\}_{i=0}^{p}$ as
\begin{equation}
  \mathcal{W}_{i}^{-1} \mathcal{L}^{(i)}\, \mathcal{W}_{i-1} = \mathcal{D} + \mathcal{E}^{(i)}.
  \label{def:LocalCocycle}
\end{equation}
We call $\mathcal{D}$ the reference transition and $\mathcal{E}^{(i)}$ the \emph{local transition perturbation}. We partition both matrices according to principal/parasitic blocks:
\begin{equation*}
  \mathcal{D} = \begin{bmatrix}
    I_{n}& \\ & D
  \end{bmatrix}, \,D:=\operatorname{diag}(\nu_2,\ldots,\nu_d)\otimes I_n,
  \quad 
  \mathcal{E}^{(i)}  =  \begin{bmatrix}
    E_{\mathrm{pr,pr}}^{(i)} & E_{\mathrm{pr,pa}}^{(i)} \\
    E_{\mathrm{pa,pr}}^{(i)} & E_{\mathrm{pa,pa}}^{(i)}
  \end{bmatrix}.
\end{equation*}
\end{definition}

The following lemma quantifies the blockwise smallness of $\mathcal{E}^{(i)}$.
\begin{lemma}[Blockwise bounds for the local transition perturbation]
  \label{lem:LocalTransitionPerturbation}
  Under \textup{(A1)}--\textup{(A4)}, there exist $h_{0}>0$ and constants $C_a,C_b$, independent of $i$, $p$ and the grid, such that for all $h\leq h_{0}$, the blocks of $\mathcal{E}^{(i)}$ satisfy
  \begin{equation*}
    \begin{aligned}
    \left\lVert E_{\mathrm{pr,pr}}^{(i)} \right\rVert + \left\lVert E_{\mathrm{pa,pr}}^{(i)} \right\rVert
    &\leq a^{(i)},\quad a^{(i)}:=C_a\,(h^{(i)})^{s+1},\\
    \left\lVert E_{\mathrm{pr,pa}}^{(i)} \right\rVert + \left\lVert E_{\mathrm{pa,pa}}^{(i)} \right\rVert
    &\leq b^{(i)},\quad b^{(i)}:=C_b\left(h^{(i)} + \sum_{\ell = 2}^{d}\left|\omega^{(i-d+\ell)}-1\right|\right).
    \end{aligned}
  \end{equation*}
\end{lemma}
\begin{proof}
By \cref{def:LocalCocycle},
\(
  \mathcal{E}^{(i)}
  =
  \mathcal{W}_{i}^{-1}
  \left(
    \mathcal{L}^{(i)}\mathcal{W}_{i-1}
    -\mathcal{W}_{i}\mathcal{D}
  \right).
\)
Recall that \(\widetilde A_j^{(i)}=\left(A_d^{(i)}\right)^{-1}A_j^{(i)},\, j=0,\ldots,d-1.\) By the structure of $\mathcal{L}^{(i)}$, only the last block row of \(\mathcal{L}^{(i)}\mathcal{W}_{i-1}-\mathcal{W}_{i}\mathcal{D}\) can be nonzero.

On the principal block column $W_{\mathrm{pr},i-1}$, the last block row equals
\[
\begin{aligned}
-\sum_{j=0}^{d-1}\widetilde A_j^{(i)}
  \Phi(t_{i-d+j})-\Phi(t_i)
&=
-\left(A_d^{(i)}\right)^{-1}
\left(
  \sum_{j=0}^{d-1}A_j^{(i)}\Phi(t_{i-d+j})
  +A_d^{(i)}\Phi(t_i)
\right) \\
&=
-\left(A_d^{(i)}\right)^{-1}
\sum_{j=0}^{d}A_j^{(i)}
  \Phi(t_{i-d+j})
=: -R_{\mathrm{pr}}^{(i)} .
\end{aligned}
\]
Since \(A_d^{(i)}=\alpha_d^{(i)}I_n-h^{(i)}\beta_d^{(i)}G(t_{i})\), \({\left(A_{d}^{(i)}\right)}^{-1}\) is uniformly bounded for sufficiently small $h$ by~\textup{(A1)}--\textup{(A2)}. Thus, using consistency assumption~\textup{(A3)}, \(\left\lVert R_{\mathrm{pr}}^{(i)} \right\rVert \leq C_{\mathrm{pr}}(h^{(i)})^{s+1}.\)

On the $\tau$-th parasitic block column $v_{\tau}$, the last block row equals
\[
\begin{aligned}
-\sum_{j=0}^{d-1}\widetilde A_j^{(i)}\nu_{\tau}^j
-\nu_{\tau}^d I_n
&=-\left(A_d^{(i)}\right)^{-1}\sum_{j=0}^{d}A_j^{(i)}\nu_{\tau}^j \\
&=-\left(A_d^{(i)}\right)^{-1}
\left[\left(\sum_{j=0}^{d}\alpha_j^{(i)}\nu_{\tau}^j\right)I_n
  -h^{(i)}\sum_{j=0}^{d}\beta_j^{(i)}G(t_{i-d+j})\nu_{\tau}^j\right] \\
&=: -R_{\mathrm{pa},\tau}^{(i)} .
\end{aligned}
\]
By (A1) and the uniform coefficient
bounds in (A2), the term proportional to $h^{(i)}$ is therefore $\mathcal{O}(h^{(i)})$. Moreover, the stability assumption \textup{(A4)} and the $C^1$-smoothness of the coefficient functions on the segment joining $\mathbf{1}$ and $\boldsymbol{\omega}^{(i)}$ together imply
\[
\left| \sum_{j=0}^{d}\alpha_j^{(i)}\nu_{\tau}^j \right| = \left|   \sum_{j=0}^{d}\left(\alpha_j(\boldsymbol{\omega}^{(i)})-\alpha_{j}(\mathbf{1})\right)\nu_{\tau}^j \right|
\leq C_\alpha \sum_{\ell = 2}^{d} \left|\omega^{(i-d+\ell)}-1\right|.
\]

Consequently, after enlarging $C_b$ if necessary,
\[
  \left\lVert R_{\mathrm{pa},\tau}^{(i)} \right\rVert
  \leq
  C_b\left(
    h^{(i)}
    +
    \sum_{\ell = 2}^{d}
    \left|\omega^{(i-d+\ell)}-1\right|
  \right).
\]

In the defect, only the last block row is nonzero: it equals $-R^{(i)}_{\mathrm{pr}}$
in the principal block column and $-R^{(i)}_{\mathrm{pa},\tau}$ in the $\tau$-th parasitic
block column. Left multiplication by $\mathcal{W}_{i}^{-1}$ and the uniform bound in
\cref{lem:MovingBasisCondition} preserve these estimates. Since $d$ is fixed, the resulting constants can be absorbed into $C_a$ and $C_b$, which gives the two stated blockwise bounds.
\end{proof}

Multiplying \cref{def:LocalCocycle} over one period and using the Floquet quasi-periodicity relation \cref{equ:MovingFrameQuasiPeriodicity} yields
\begin{equation}
  \mathcal{W}_{0}^{-1} \mathcal{F}_{p,0}\, \mathcal{W}_{0}
  =
  \mathcal{B}\,\left(\mathcal{D}+\mathcal{E}^{(p)}\right)\cdots\left(\mathcal{D}+\mathcal{E}^{(1)}\right).
  \label{equ:OnePeriodCocycle}
\end{equation}
In this basis, the one-period propagator $\mathcal{F}_{p,0}$ is therefore a product of almost block-diagonal factors.

\subsection{Stepwise Block Triangularization via a Periodic Riccati Equation}
\label{subsec:Triangularization}

We seek a stepwise block upper triangularization, compatible over one period, that isolates the perturbed principal and parasitic dynamics:
\begin{equation}
  \left(\mathcal{T}_{i}\right)^{-1}
  \left(\mathcal{D}+\mathcal{E}^{(i)}\right)
  \mathcal{T}_{i-1}
  =
  \begin{bmatrix}
    P_{i} & * \\
    0 & H_{i}
  \end{bmatrix},\quad i = 1,\ldots,p,\quad
  \mathcal{T}_{0}\mathcal{B}\left(\mathcal{T}_{p}\right)^{-1} = \mathcal{B},
  \label{equ:LocalTriangularization}
\end{equation}
where
\begin{equation}
  \mathcal{T}_{i}
  :=
  \begin{bmatrix}
    I_n & 0 \\
    \Gamma_{i} & I_{(d-1)n}
  \end{bmatrix},
  \qquad
  \left(\mathcal{T}_{i}\right)^{-1}
  =
  \begin{bmatrix}
    I_n & 0 \\
    -\Gamma_{i} & I_{(d-1)n}
  \end{bmatrix}.
  \label{def:RiccatiTransformation}
\end{equation}
The diagonal blocks in \cref{equ:LocalTriangularization} are defined by
\begin{equation}
  P_{i}
  :=
  I_n+E_{\mathrm{pr,pr}}^{(i)}
  +E_{\mathrm{pr,pa}}^{(i)}\Gamma_{i-1},
  \qquad
  H_{i}
  :=
  D+E_{\mathrm{pa,pa}}^{(i)}
  -\Gamma_{i}E_{\mathrm{pr,pa}}^{(i)}.
  \label{def:PtH}
\end{equation}
The lower-left block in \cref{equ:LocalTriangularization} vanishes if and only if
\begin{equation}
  E_{\mathrm{pa,pr}}^{(i)}
  +\left(D+E_{\mathrm{pa,pa}}^{(i)}\right)\Gamma_{i-1}
  -\Gamma_{i}
   \left(
     I_n+E_{\mathrm{pr,pr}}^{(i)}
     +E_{\mathrm{pr,pa}}^{(i)}\Gamma_{i-1}
   \right)
  =0.
  \label{equ:DiscreteRiccatiEquation}
\end{equation}
Periodicity of the transformation additionally imposes the Floquet boundary condition
\begin{equation}
  \Gamma_{p} = \Gamma_{0}\Phi(T,0).
  \label{equ:FloquetBoundaryCondition}
\end{equation}
Indeed, $\mathcal{T}_{0}\mathcal{B}(\mathcal{T}_{p})^{-1}=\mathcal{B}$ means $\mathcal{B}\mathcal{T}_{p}=\mathcal{T}_{0}\mathcal{B}$, which by \cref{def:RiccatiTransformation} is exactly
\cref{equ:FloquetBoundaryCondition}.

When each $P_{i}$ is invertible, the recursion \cref{equ:DiscreteRiccatiEquation} advances $\Gamma_{i-1}$ to $\Gamma_{i}$. Imposing the boundary condition \cref{equ:FloquetBoundaryCondition} then gives a Riccati boundary-value problem. The next lemma shows that this Riccati boundary-value problem has a unique small solution satisfying the stated bound, uniformly over the admissible grids. The bound on $\Gamma_{i}$ is then used to estimate the spectrum of $\mathcal{F}_{p,0}$ via \cref{equ:LocalTriangularization}.

\begin{lemma}[Small solution of the periodic Riccati boundary-value problem]
  \label{lem:PeriodicGraph}
  Let \textup{(A1)}--\textup{(A4)} hold. For all sufficiently small
  maximum stepsize $h$, there exists a sequence
  $\Gamma_{i}\in\mathbb{C}^{(d-1)n\times n}$, $i=0,\ldots,p$, satisfying
  \begin{equation}
    \Gamma_{i}
    =
    \left[
      E_{\mathrm{pa,pr}}^{(i)}
      +
      \left(D+E_{\mathrm{pa,pa}}^{(i)}\right)\Gamma_{i-1}
    \right]
    \left[
      I_n
      +E_{\mathrm{pr,pr}}^{(i)}
      +E_{\mathrm{pr,pa}}^{(i)}\Gamma_{i-1}
    \right]^{-1},
    \qquad i=1,\ldots,p,
    \label{equ:GraphRecursion}
  \end{equation}
  together with the Floquet boundary condition
  \cref{equ:FloquetBoundaryCondition}. Moreover, there is a constant
  $C_\Gamma$, independent of the grid, such that
  \begin{equation}
    \max_{0\leq i\leq p}\left\lVert \Gamma_{i} \right\rVert
    \leq C_\Gamma h^{s+1}.
    \label{equ:GraphBound}
  \end{equation}
  The solution is unique among sequences satisfying
  \cref{equ:GraphBound} with this constant $C_\Gamma$.
\end{lemma}

For the proof of this lemma and the subsequent product estimates, use the following shorthand for the bounds in
\cref{lem:LocalTransitionPerturbation}. For $1\leq i\leq p$, set
\[
a_\infty:=C_a h^{s+1}, \qquad A:=C_aT,\qquad B:=C_b(T+dC_T).
\]
Thus, 
\begin{equation}
  \max_{1\leq i\leq p}a^{(i)} \leq a_\infty,\qquad
\sum_{i=1}^{p}a^{(i)} \leq Ah^{s} ,\qquad
\sum_{i=1}^{p}b^{(i)}\leq B.
\label{equ:WeightBounds}
\end{equation}
In particular, every interval sum of the weights $b^{(i)}$ is bounded by $B$.

\begin{proof}
The proof applies the Banach contraction principle \cite[Theorem~1.1, p.~3]{Pata2019FixedPoint}. The recurrence for $\Gamma_{i}$ induces a map from $\Gamma_{0}$ to $\Gamma_{p}$, and the Floquet boundary condition turns it into a self-map of the admissible set of initial values. A fixed point of this map solves the periodic Riccati boundary-value problem. The following claim isolates the stability estimates needed to prove that this self-map is contractive; its proof is deferred until after the fixed-point argument.
\begin{claim}[Stability of the Riccati recursion]
  \label{clm:GraphStability}
  Recall \(r:=\left\lVert D \right\rVert=\max_{2\leq\tau\leq d}|\nu_{\tau}|<1.\) There exist constants $q\in(r,1)$ and $C_{\mathrm s}\ge1$, independent of the grid, such that the following holds. For every fixed $R>0$, there exists $h_{R}>0$ such that, for all $h\leq h_{R}$, the recursion \cref{equ:GraphRecursion} is well defined whenever the preceding iterates lie in
  \(
    \mathbb{B}_{Ra_{\infty}}:=\{\Gamma:\left\lVert \Gamma \right\rVert\leq Ra_{\infty}\}.
  \)
  Moreover, every orbit segment contained in this ball satisfies
  \begin{equation}
    \left\lVert \Gamma_{i} \right\rVert
    \leq
    C_{\mathrm s}q^{\,i-j}\left\lVert \Gamma_{j} \right\rVert
    +
    C_{\mathrm s}
    \sum_{k=j+1}^{i}q^{\,i-k}a^{(k)},
    \qquad j<i,
    \label{equ:GraphPropagation}
  \end{equation}
  and any two such orbit segments satisfy
  \begin{equation}
    \left\lVert \Gamma_{i}-\Gamma^{\prime}_{i} \right\rVert
    \leq
    C_{\mathrm s}q^{\,i-j}
    \left\lVert \Gamma_{j}-\Gamma^{\prime}_{j} \right\rVert.
    \label{equ:GraphDifferencePropagation}
  \end{equation}
\end{claim}

Set $\kappa_\Phi:=\left\lVert \Phi(T,0)^{-1} \right\rVert$ and \(\kappa_{q}:=C_{\mathrm s}/(1-q)\), and define \(R_{0}:=2\kappa_\Phi \kappa_{q}\) and \(R_1:=C_{\mathrm s}R_{0}+\kappa_{q}\). The ball of radius $R_{0}a_\infty$ is the domain of the fixed-point map, whereas the ball of radius $R_1a_\infty$ contains the one-period trajectories generated from it.

Apply \cref{clm:GraphStability} with $R=R_1$. Since $p\ge T/h\to\infty$, we have $q^p\to0$; hence, for sufficiently small $h$,
\begin{equation}
  \kappa_\Phi C_{\mathrm s}q^p\leq\frac{1}{2}.
  \label{equ:PeriodContractionSmallness}
\end{equation}

\medskip
\noindent
\emph{Step 1: trajectories starting in the fixed-point ball $\mathbb{B}_{R_{0}a_{\infty}}$ remain in the larger stability ball $\mathbb{B}_{R_{1}a_{\infty}}$ .}
Suppose that $\left\lVert \Gamma_{0} \right\rVert\leq R_{0}a_\infty$.  As long as the preceding iterates remain in $\mathbb{B}_{R_1a_{\infty}}$,
\cref{equ:GraphPropagation} gives
\begin{equation}
  \begin{aligned}
    \left\lVert \Gamma_{i} \right\rVert
    &\leq
    C_{\mathrm s}q^i\left\lVert \Gamma_{0} \right\rVert
    +
    C_{\mathrm s}\sum_{k=1}^{i}q^{\,i-k}a^{(k)}
    \leq
    C_{\mathrm s}R_{0}a_\infty
    +
    C_{\mathrm s}\frac{1-q^i}{1-q}a_\infty
    \\
    &\leq
    \left(C_{\mathrm s}R_{0}+\kappa_{q}\right)a_\infty
    =
    R_1a_\infty.
  \end{aligned}
  \label{equ:OrbitBound}
\end{equation}
Since $\Gamma_{0}\in\mathbb{B}_{R_1a_{\infty}}$, the preceding estimate and \cref{clm:GraphStability} imply inductively that the recursion is well defined and \(\Gamma_{i}\in\mathbb{B}_{R_{1}a_{\infty}},\, i=0,\ldots,p.\)

\medskip
\noindent
\emph{Step 2: the Floquet-adjusted one-period map returns to the fixed-point ball $\mathbb{B}_{R_{0}a_{\infty}}$.}
For $\Gamma_{0}\in\mathbb{B}_{R_{0}a_{\infty}}$, define
\(
  \mathfrak{G}_{h}(\Gamma_{0}):=\Gamma_{p}(\Gamma_{0})\Phi(T,0)^{-1}.
\)
The fixed points of $\mathfrak{G}_{h}$ are exactly the solutions of \cref{equ:GraphRecursion,equ:FloquetBoundaryCondition}. By \cref{equ:GraphPropagation} and \cref{equ:PeriodContractionSmallness},
\begin{align*}
  \left\lVert \mathfrak{G}_{h}(\Gamma_{0}) \right\rVert
  &\leq
  \kappa_\Phi
  \left(C_{\mathrm s}q^pR_{0}+\kappa_{q}\right)a_\infty
  \leq
  \left(\frac12R_{0}+\kappa_\Phi \kappa_{q}\right)a_\infty
  =
  R_{0}a_\infty.
\end{align*}
Thus, $\mathfrak{G}_{h}$ maps $\mathbb{B}_{R_{0}a_{\infty}}$ into itself.

\emph{Step 3: the iteration is a contraction.}
For $\Gamma_{0},\Gamma^{\prime}_{0}\in\mathbb{B}_{R_{0}a_{\infty}}$, their trajectories lie in $\mathbb{B}_{R_{1}a_{\infty}}$ by \cref{equ:OrbitBound}; hence \cref{equ:GraphDifferencePropagation} in \cref{clm:GraphStability} gives
\[
  \|\mathfrak{G}_{h}(\Gamma_{0})
    -\mathfrak{G}_{h}(\Gamma^{\prime}_{0})\|
  \leq
  \kappa_\Phi C_{\mathrm s}q^p
  \left\lVert \Gamma_{0}-\Gamma^{\prime}_{0} \right\rVert
  \leq
  \frac12\left\lVert \Gamma_{0}-\Gamma^{\prime}_{0} \right\rVert.
\]
Therefore, $\mathfrak{G}_{h}$ is a contraction on $\mathbb{B}_{R_{0}a_{\infty}}$.  The Banach contraction principle \cite[Theorem~1.1, p.~3]{Pata2019FixedPoint} yields a solution of \cref{equ:GraphRecursion,equ:FloquetBoundaryCondition}, and \cref{equ:OrbitBound} gives
\[
  \max_{0\leq i\leq p}\left\lVert \Gamma_{i} \right\rVert
  \leq R_1a_\infty
  =C_aR_1h^{s+1}.
\]
Thus, \cref{equ:GraphBound} holds with $C_\Gamma:=C_aR_1$.

Finally, suppose that $\Gamma_{i}$ and $\Gamma^{\prime}_{i}$ are two solutions of \cref{equ:GraphRecursion,equ:FloquetBoundaryCondition} satisfying \cref{equ:GraphBound} with this $C_\Gamma$. Both trajectories lie in $\mathbb{B}_{R_{1}a_{\infty}}$, so \cref{equ:GraphDifferencePropagation} and the Floquet boundary condition give
\[
  \left\lVert \Gamma_{0}-\Gamma^{\prime}_{0} \right\rVert
  \leq
  \kappa_\Phi C_{\mathrm s}q^p
  \left\lVert \Gamma_{0}-\Gamma^{\prime}_{0} \right\rVert
  \leq
  \frac12\left\lVert \Gamma_{0}-\Gamma^{\prime}_{0} \right\rVert.
\]
Hence $\Gamma_{0}=\Gamma^{\prime}_{0}$, and then \cref{equ:GraphRecursion} implies $\Gamma_{i}=\Gamma^{\prime}_{i}$ for every $i$, which shows the uniqueness of the solution.
\end{proof}

\begin{proof}[Proof of \cref{clm:GraphStability}]
Set $q:=(1+r)/2$. For fixed $R>0$, write
\[
\xi^{(i)}:=a^{(i)}+Ra_\infty b^{(i)},\qquad
N_{i}(\Gamma):=
I_n+E_{\mathrm{pr,pr}}^{(i)}
+E_{\mathrm{pr,pa}}^{(i)}\Gamma .
\]
Since
$\sum_{i=1}^{p}\xi^{(i)}\leq Ah^{s}+Ra_\infty B=\mathcal{O}(h^s)$, we may assume, after decreasing $h$ if necessary, that $\sum_i\xi^{(i)}\le1/4$. Hence, for $\Gamma\in\mathbb{B}_{Ra_{\infty}}$, $\left\lVert N_{i}(\Gamma)-I_n \right\rVert\leq\xi^{(i)}\le1/4$ and $\left\lVert N_{i}(\Gamma)^{-1} \right\rVert\leq(1-\xi^{(i)})^{-1} \leq e^{2\xi^{(i)}}$.

Moreover, $r+b^{(k)}\leq q(1+b^{(k)}/q)$. Since every interval sum of the nonnegative weights $b^{(i)}$ is bounded by $ B$ \cref{equ:WeightBounds}, we obtain, for each fixed $c>0$,
\begin{equation}
\prod_{k=\ell}^{i}
e^{c\xi^{(k)}}(r+b^{(k)})
\leq
e^{\frac{c}{4}} q^{\,i-\ell+1}e^{\sum_{k = \ell}^{i}\frac{b^{(k)}}{q}}\leq e^{\frac{c}{4}+\frac{ B}{q}}q^{\,i-\ell+1}
=:C_cq^{\,i-\ell+1}.
\label{equ:UniformGraphProductBound}
\end{equation}

Let $\mathfrak{f}_{i}$ denote the right-hand side of \cref{equ:GraphRecursion}. The preceding estimates give \[\left\lVert \mathfrak{f}_{i}(\Gamma) \right\rVert \leq e^{2\xi^{(i)}}a^{(i)} +e^{2\xi^{(i)}}(r+b^{(i)})\left\lVert \Gamma \right\rVert.\]
Iterating this inequality from $j$ to $i$ and applying \cref{equ:UniformGraphProductBound} with $c=2$ yields
\[
\begin{aligned}
\left\lVert \Gamma_{i} \right\rVert
\leq&\prod_{k=j+1}^{i}e^{2\xi^{(k)}}(r+b^{(k)})\left\lVert \Gamma_{j} \right\rVert+\sum_{k=j+1}^{i}e^{2\xi^{(k)}}a^{(k)}\prod_{\ell=k+1}^{i}e^{2\xi^{(\ell)}}(r+b^{(\ell)})
\\
\leq&
C_1q^{\,i-j}\left\lVert \Gamma_{j} \right\rVert
+
C_1\sum_{k=j+1}^{i}q^{\,i-k}a^{(k)},
\end{aligned}
\]
where $C_1$ is independent of the grid.

For two points $\Gamma,\Gamma^{\prime}\in\mathbb{B}_{Ra_{\infty}}$, we have
\[
\begin{aligned}
\mathfrak{f}_{i}(\Gamma)-\mathfrak{f}_{i}(\Gamma^{\prime})
&=
\left(D+E_{\mathrm{pa,pa}}^{(i)}\right)
(\Gamma-\Gamma^{\prime})N_{i}(\Gamma)^{-1}\\
&\quad+
\left[
E_{\mathrm{pa,pr}}^{(i)}
+\left(D+E_{\mathrm{pa,pa}}^{(i)}\right)\Gamma^{\prime}
\right]
\left[
N_{i}(\Gamma)^{-1}
-N_{i}(\Gamma^{\prime})^{-1}
\right].
\end{aligned}
\]
By \cref{lem:LocalTransitionPerturbation} and the definition of
$\mathbb{B}_{Ra_{\infty}}$,
\(
\left\lVert D+E_{\mathrm{pa,pa}}^{(i)}\right\rVert
\leq r+b^{(i)},
\,
\left\lVert E_{\mathrm{pa,pr}}^{(i)}\right\rVert
\leq a^{(i)},
\,
\left\lVert \Gamma^{\prime}\right\rVert
\leq Ra_\infty.
\)
Moreover, the resolvent identity, together with the bounds for
$N_{i}(\Gamma)^{-1}$ and $N_{i}(\Gamma^{\prime})^{-1}$, gives
\[
\begin{aligned}
\left\lVert
N_{i}(\Gamma)^{-1}
-N_{i}(\Gamma^{\prime})^{-1}
\right\rVert &\leq
\left\lVert N_{i}(\Gamma)^{-1}\right\rVert
\left\lVert
N_{i}(\Gamma)-N_{i}(\Gamma^{\prime})
\right\rVert
\left\lVert N_{i}(\Gamma^{\prime})^{-1}\right\rVert\\
&\leq
b^{(i)}e^{4\xi^{(i)}}
\left\lVert\Gamma-\Gamma^{\prime}\right\rVert.
\end{aligned}
\]
Therefore,
\[
\begin{aligned}
\left\lVert
\mathfrak{f}_{i}(\Gamma)
-\mathfrak{f}_{i}(\Gamma^{\prime})
\right\rVert
&\leq
(r+b^{(i)})e^{2\xi^{(i)}}
\left\lVert\Gamma-\Gamma^{\prime}\right\rVert
\\
&\quad+
b^{(i)}
\left(
a^{(i)}
+(r+b^{(i)})Ra_\infty
\right)
e^{4\xi^{(i)}}
\left\lVert\Gamma-\Gamma^{\prime}\right\rVert.
\end{aligned}
\]
Since \((r+b^{(i)})\xi^{(i)}- b^{(i)}\left(a^{(i)}+(r+b^{(i)})Ra_\infty\right) =ra^{(i)}\ge0,\) we have 
\[b^{(i)} \left(a^{(i)}+(r+b^{(i)})Ra_\infty\right)\leq (r+b^{(i)})\xi^{(i)}.\]
Moreover, because \(0\leq\xi^{(i)}\le1/4\), \(e^{2\xi^{(i)}}\leq e^{1/2}<2 \) and hence \(1+\xi^{(i)}e^{2\xi^{(i)}} \leq 1+2\xi^{(i)} \leq e^{2\xi^{(i)}}.\) Consequently,
\[
\begin{aligned}
\left\lVert
\mathfrak{f}_{i}(\Gamma)
-\mathfrak{f}_{i}(\Gamma^{\prime})
\right\rVert
&\leq
(r+b^{(i)})e^{2\xi^{(i)}}
\left(
1+\xi^{(i)}e^{2\xi^{(i)}}
\right)
\left\lVert\Gamma-\Gamma^{\prime}\right\rVert \\
&\leq
(r+b^{(i)})e^{4\xi^{(i)}}
\left\lVert\Gamma-\Gamma^{\prime}\right\rVert.
\end{aligned}
\]
Iterating this one-step estimate from \(j\) to \(i\) and applying \cref{equ:UniformGraphProductBound} with \(c=4\) yields
\[
\left\lVert \Gamma_{i}-\Gamma^{\prime}_{i} \right\rVert \leq
\prod_{k=j+1}^{i} \left[(r+b^{(k)})e^{4\xi^{(k)}} \right]\left\lVert
\Gamma_{j}-\Gamma^{\prime}_{j} \right\rVert \leq C_2q^{\,i-j}
\left\lVert \Gamma_{j}-\Gamma^{\prime}_{j} \right\rVert,
\]
where \(C_2\) is independent of the grid. Taking
$C_{\mathrm s}:=\max\{1,C_1,C_2\}$ closes the proof.
\end{proof}

\subsection{Proof of the Main Convergence Theorem}
\label{subsec:ProofMain}
For the subspace estimate in part~(3), we use the following specialization of Stewart's invariant-subspace perturbation theorem \cite[Theorem~4.1]{Stewart1971ISbounds}; see also \cite{Betcke2011PerturbationNew}.
\begin{theorem}[Invariant Subspace Perturbation Theorem \cite{Stewart1971ISbounds,Betcke2011PerturbationNew}]
  \label{thm:Perturbation}
  Let
  \[
    \Lambda=\operatorname{diag}(\Lambda_m,\Lambda_{n-m}), \qquad E=
    \begin{bmatrix}
      E_{11}&E_{12}\\
      E_{21}&E_{22}
    \end{bmatrix},
  \]
  and let \(\gamma_m\) be defined by \cref{def:ClusterSep}. Set \(\delta:=\gamma_m-\left(\lVert E_{11}\rVert_F+\lVert E_{22}\rVert_F\right).\)   If \(\delta>0\), \(\lVert E_{21}\rVert_F\lVert E_{12}\rVert_F<\frac{\delta^2}{4},\) then \(\Lambda+E\) has an invariant subspace \(\operatorname{range}(Y_m)\), where
  \[
    Y_m=
    \begin{bmatrix}
      I_m\\K_m
    \end{bmatrix},
    \qquad
    \lVert K_m\rVert_F
    \leq \frac{2\lVert E_{21}\rVert_F}{\delta}.
  \]
\end{theorem}

\begin{proof}[Proof of \cref{thm:MainConvergence}]
Let \(\{\Gamma_i\}_{i=0}^{p}\) denote the solution sequence constructed in \cref{equ:GraphRecursion}, which satisfies the Floquet boundary condition. By \cref{equ:LocalTriangularization}, the transformations $\mathcal{T}_i$ block triangularize every local transition; since $\mathcal{B}\mathcal{T}_p=\mathcal{T}_0\mathcal{B}$, they are compatible over one period.

Define \( P_{p,0}:=P_{p}\cdots P_{1},\) \(H_{p,0}:=H_{p}\cdots H_{1},\) and introduce the fixed Floquet coordinate transformation
\begin{equation*}
\mathcal{Q}:=\operatorname{diag}\left(U(0),I_{(d-1)n}\right), \qquad \tilde{P}_{p,0}:=U(0)^{-1}P_{p,0}U(0).
\end{equation*}
Multiplying \cref{equ:LocalTriangularization} from \(i=1\) to \(p\),
using the one-period representation \cref{equ:OnePeriodCocycle}, and then using \(U(0)^{-1}\Phi(T,0)U(0)=\Lambda\), we obtain the exact similarity
\begin{equation}
\left(\mathcal{W}_{0}\mathcal{T}_{0}\mathcal{Q}\right)^{-1}
\mathcal{F}_{p,0}
\left(\mathcal{W}_{0}\mathcal{T}_{0}\mathcal{Q}\right)
=
\begin{bmatrix}
\Lambda\tilde{P}_{p,0} & *\\
0 & H_{p,0}
\end{bmatrix}.
\label{equ:OnePeriodSimilarity}
\end{equation}
Consequently,
\begin{equation}
\sigma(\mathcal{F}_{p,0}) = \sigma(\Lambda\tilde{P}_{p,0}) \cup \sigma(H_{p,0}).
\label{equ:PrincipalParasiticSpectralSplit}
\end{equation}

We next estimate the two diagonal products.  From \cref{def:PtH,lem:LocalTransitionPerturbation,equ:GraphBound},
\[
\begin{aligned}
\eta_{h} &:=\sum_{i=1}^{p}\left\lVert P_{i}-I_n\right\rVert  \leq \sum_{i=1}^{p}a^{(i)} + \left( \max_{0\leq i\leq p}\left\lVert\Gamma_{i}\right\rVert \right)
\sum_{i=1}^{p}b^{(i)} \\
&\leq Ah^{s}+C_\Gamma h^{s+1} B =\mathcal{O}(h^s).
\end{aligned}
\]
Using the ordered-product estimate
\(
\left\lVert\prod_i(I+Z_i)-I\right\rVert \leq \exp\left(\sum_i\left\lVert Z_i\right\rVert\right)-1
\)
gives
\[
\left\lVert P_{p,0}-I_n\right\rVert \leq e^{\eta_{h}}-1 =\mathcal{O}(h^s).
\]
Since \(U(0)\) is fixed, it follows that
\begin{equation}
S_{h} :=\Lambda\tilde{P}_{p,0} =\Lambda+E_{\Lambda},
\qquad E_{\Lambda} :=\Lambda\left(\tilde{P}_{p,0}-I_n\right),
\qquad \varepsilon_{h}:=\left\lVert E_{\Lambda}\right\rVert_{\mathrm{F}} \leq C_{E} h^s,
\label{equ:PrincipalPerturbationBound}
\end{equation}
where \(C_{E}\) is independent of the grid.

For the parasitic product, \cref{def:PtH} and \cref{lem:LocalTransitionPerturbation,equ:GraphBound} imply
\[
\left\lVert H_{i}\right\rVert \leq r+w^{(i)}, \qquad w^{(i)}:=\left(1+C_\Gamma h^{s+1}\right)b^{(i)}.
\]
After decreasing \(h\) if necessary, \(\sum_{i=1}^{p}w^{(i)}\le2 B.\) If \(r>0\), then
\begin{equation}
\left\lVert H_{p,0}\right\rVert \leq \prod_{i=1}^{p}\left(r+w^{(i)}\right) \leq r^p \exp\left(\frac{1}{r}\sum_{i=1}^{p}w^{(i)}\right) \leq C r^p, \quad C:=e^{2 B/r}.
\label{equ:ParasiticPositiveRootBound}
\end{equation}

We can now prove the three assertions.

\emph{Part~(1).}
The matrix \(\Lambda\) is diagonal and hence normal.  By
\cref{equ:PrincipalPerturbationBound} and the Bauer--Fike theorem for normal matrices \cite{bauer1960NormsExclusion}, every eigenvalue of \(S_{h}\) lies within \(\mathcal{O}(h^s)\) of \(\sigma(\Lambda)\). For sufficiently small \(h\), the neighborhoods belonging to distinct Floquet multipliers are disjoint, and continuity of the characteristic roots, counted with algebraic multiplicity, preserves the number of eigenvalues in each neighborhood.  Therefore, after suitable ordering,
\[
\widehat{\lambda}_{k}=\lambda_{k}+\mathcal{O}(h^s),
\qquad k=1,\ldots,n.
\]

\emph{Part~(2).}
By \cref{equ:PrincipalParasiticSpectralSplit}, the parasitic periodic eigenvalues are the eigenvalues of \(H_{p,0}\).  Their moduli are bounded by \(\lVert H_{p,0}\rVert\), so \cref{equ:ParasiticPositiveRootBound} proves the estimate in part~(2).  Since every Floquet multiplier is nonzero, parts~(1)--(2) also show that the principal and parasitic spectral groups are disjoint for all sufficiently small \(h\).

\emph{Part~(3).}
Partition \(E_{\Lambda}=[E_{ij}]_{i,j=1}^2\) conformably with \(\Lambda=\operatorname{diag}(\Lambda_m,\Lambda_{n-m})\). Since \(\lVert E_{ij}\rVert_F\leq\varepsilon_{h}\) and \(\varepsilon_{h}=\mathcal{O}(h^s)\), for all sufficiently small \(h\),
\[
\begin{aligned}
\delta_{h} &:=\gamma_m-\left(\lVert E_{11}\rVert_F+\lVert E_{22}\rVert_F\right)
\ge \gamma_m-2\varepsilon_{h}
>\frac{\gamma_m}{2},\quad
\lVert E_{21}\rVert_F\lVert E_{12}\rVert_F &\leq \varepsilon_{h}^2<\frac{\delta_{h}^2}{4}.
\end{aligned}
\]
Applying \cref{thm:Perturbation} with \(E=E_{\Lambda}\) gives an invariant subspace \(\operatorname{range}(Y_m)\) of \(S_{h}=\Lambda+E_{\Lambda}\), where
\[
Y_m=\begin{bmatrix}I_m\\K_m\end{bmatrix},
\qquad\lVert K_m\rVert_F
\leq \frac{2\lVert E_{21}\rVert_F}{\delta_{h}}
\leq \frac{4\varepsilon_{h}}{\gamma_m}
=\mathcal{O}\left(\frac{h^s}{\gamma_m}\right).
\]
The restriction of \(S_{h}\) to this invariant subspace is represented by $\Lambda_m+E_{11}+E_{12}K_m,$ which is an \(\mathcal O(h^s)\) perturbation of \(\Lambda_m\). Hence, by part~(1) and \(\gamma_m>0\), the subspace is the spectral subspace associated with the \(m\) eigenvalues converging to \(\sigma(\Lambda_m)\).

The block triangularization in \cref{equ:OnePeriodSimilarity} therefore lifts this subspace to
\[
\begin{aligned}
\widehat{\mathcal{U}}_m
&=\operatorname{range}\left(
\mathcal{W}_{0}\mathcal{Q}
\begin{bmatrix}
I_m\\K_m\\\Gamma_{0}U(0)Y_m
\end{bmatrix}
\right),\quad
\mathcal{U}_m
&=\operatorname{range}\left(
\mathcal{W}_{0}\mathcal{Q}
\begin{bmatrix}
I_m\\0\\0
\end{bmatrix}
\right).
\end{aligned}
\]
The bound for \(K_m\) implies \(\lVert Y_m\rVert=\left\lVert\begin{bmatrix}
  I_{m}\\K_{m}
\end{bmatrix} \right\rVert =\mathcal{O}(1)\). Together with \cref{equ:GraphBound}, it follows that
\[
\left\lVert
\begin{bmatrix}
K_m\\\Gamma_{0}U(0)Y_m
\end{bmatrix}
\right\rVert
=\mathcal{O}\left(\frac{h^s}{\gamma_m}+h^{s+1}\right)
=\mathcal{O}\left(\frac{h^s}{\gamma_m}\right),
\]
where the last equality uses the fact that \(\gamma_m\) is fixed and positive. Finally, \cref{lem:MovingBasisCondition} and the fixed invertibility of \(\mathcal{Q}\) show that \(\mathcal{W}_{0}\mathcal{Q}\) is uniformly conditioned. Then, using the bound for canonical angles between subspaces \cite[Eq.(2.8)]{Stewart1971ISbounds} and applying the uniformly conditioning of \(\mathcal{W}_{0}\mathcal{Q}\) gives
\[
\left\lVert \sin \Theta\left(
\widehat{\mathcal{U}}_m,\mathcal{U}_m
\right) \right\rVert =\mathcal{O}\left(\frac{h^s}{\gamma_m}\right).
\]
Thus, the part~(3) has been proved.
\end{proof}

\section{An Efficient Solver of the pPEP}
\label{sec:Algorithm}
To compute the dominant periodic eigenvalues of the pPEP \cref{equ:pPEPDefinition}, one may apply the periodic Krylov--Schur method (\texttt{pKS}) \cite{kressner2006PeriodicKS} to its linearization \cref{equ:EVPLinearizedPeriodic}. The \texttt{pKS} constructs a periodic Arnoldi process
\begin{equation}
  \mathcal{L}^{(i)}V^{(i-1)} = V^{(i)}T^{(i)},\quad i = 1,\ldots,p,
  \label{equ:periodicKS}
\end{equation}
where $V^{(1)},\ldots,V^{(p-1)} \in \mathbb{R}^{dn\times k}$ and $V^{(p)} \in \mathbb{R}^{dn\times (k+1)}$ have orthonormal columns, with $V^{(0)}$ denoting the first $k$ active columns of $V^{(p)}$. Accordingly, $T^{(i)}\in\mathbb R^{k\times k}$ for $1\le i\le p-1$, and $T^{(p)}\in\mathbb R^{(k+1)\times k}$. These conforming factors form a periodic Hessenberg decomposition, from which Ritz values are extracted. As a general periodic eigenvalue solver, \texttt{pKS} exploits the periodic structure and is numerically stable \cite{kressner2006PeriodicKS}.

However, the linearization of the pPEP increases the problem dimension from $n$ to $dn$. The companion structure can be exploited to reduce the memory usage.

Motivated by the Two-level Orthogonal Arnoldi (TOAR) procedure \cite{lu2016TOAR}, we propose \texttt{pTOAR}, which compresses the periodic Arnoldi bases $V^{(i)}$ in \texttt{pKS} via a two-level orthogonal factorization. The proposed two-level factorization reformulates the periodic Arnoldi expansion while preserving the periodic Hessenberg relation used for Ritz extraction.

\subsection{Periodic Two-level Orthogonal Arnoldi Process (\texttt{pTOAR})}
The main idea of \texttt{pTOAR} is to represent $V^{(i)}$ implicitly in compact form when the matrices $\mathcal{L}^{(i)}$ have the companion structure in \cref{def:L}. With $\widetilde A_j^{(i)}:=\left(A_d^{(i)}\right)^{-1}A_j^{(i)}$, we partition $V^{(i)}$ according to this block structure:
\begin{equation}\begin{bmatrix}
    &I_{n}&&\\
    &&\ddots&\\
    &&&I_{n}\\
    -\widetilde A_{0}^{(i)}&-\widetilde A_{1}^{(i)}&\cdots&-\widetilde A_{d-1}^{(i)}\\
  \end{bmatrix}\begin{bmatrix}
    V_1^{(i-1)}\\ \vdots\\V_{d-1}^{(i-1)} \\ V_d^{(i-1)}
  \end{bmatrix}=\begin{bmatrix}
    V_1^{(i)}\\ \vdots\\V_{d-1}^{(i)} \\ V_d^{(i)}
  \end{bmatrix}T^{(i)}.
\label{equ:pTOARpre}
\end{equation}
From \cref{equ:pTOARpre}, we obtain \[V_j^{(i)}T^{(i)}=V_{j+1}^{(i-1)}, \qquad j=1,\ldots,d-1,\quad i=1,\ldots,p,
\]
with periodic wrapping in the superscript. This suggests introducing, for each phase, a shared orthonormal basis $Q^{(i)}$ such that
\[
\operatorname{span}\left(Q^{(i)}\right) = \operatorname{span}\left(V_d^{(i)},V_{d-1}^{(i+1)},\ldots,V_{1}^{(i+d-1)}\right).
\]
With this shared basis, each Arnoldi basis $V^{(i)}$ can be compactly expressed:
$$
V^{(i)} = 
\begin{bmatrix}
  Q^{(i-d+1)} & & \\
  & \ddots & \\
  & & Q^{(i)}
\end{bmatrix}
\begin{bmatrix}
  U_1^{(i)} \\
  \vdots \\
  U_d^{(i)}
\end{bmatrix},
$$
where the $U_j^{(i)}$ are local representation matrices of conforming sizes. Because the block-diagonal $Q$ factor has orthonormal columns and $V^{(i)}$ is an orthonormal Arnoldi basis, the compact representation implies that $\operatorname{col}(U_1^{(i)},\ldots,U_d^{(i)})$ also has orthonormal columns. For notational simplicity, we suppress the phase-dependent boundary dimensions and interpret all compact factors with conforming sizes. 

We now focus on the periodic Arnoldi iteration. A \texttt{pKS} expansion performs a full cycle over $i=1,\ldots,p$. The companion structure implies that only the last block of the new Arnoldi vector must be computed explicitly.
\begin{equation}
  \mathcal{L}^{(i)}v^{(i-1)} = \left[
{v_{2}^{(i-1)}}^{\mathrm{H}} ,\cdots, {v_{d}^{(i-1)}}^{\mathrm{H}} , {v_{d}^{(i)}}^{\mathrm{H}}
  \right]^{\mathrm{H}}, \, v_{d}^{(i)}=-{A_d^{(i)}}^{-1}\left(\sum_{j=0}^{d-1} A_j^{(i)} v_{j+1}^{(i-1)}\right),
\label{equ:pKSv}
\end{equation}
where $v^{(i-1)}$ is the last column of $V^{(i-1)}$, and $v_{j}^{(i-1)}$ denotes its $j$-th block, equivalently the last column of $V_{j}^{(i-1)}$ in \cref{equ:pTOARpre}. 

In pTOAR, the preceding blocks $v_j^{(i-1)}$ are reconstructed from the compressed factors $Q^{(i-d+j-1)}$ and $U_j^{(i-1)}$, after which \cref{equ:pKSv} gives the new last block. The first-level orthogonalization expands the shared basis $Q^{(i)}$, whereas the second level updates the coefficient matrices $U_j^{(i)}$ so that the implicitly represented basis $V^{(i)}$ remains orthonormal. 

\Cref{alg:pTOAR} outlines one full periodic Arnoldi expansion sweep of \texttt{pTOAR} in compact form. Deflation and breakdown handling are outside the present scope.
\begin{algorithm}[ht]
\caption{\texttt{pTOAR}: periodic Arnoldi process with TOAR-style compression}
\label{alg:pTOAR}
\begin{algorithmic}[1]
\REQUIRE
$\{A^{(i)}_j\}_{j=0}^d$ for $i=1,\ldots,p$ and current conforming compact factors $\{Q^{(i)}\}$, $\{U_j^{(i)}\}$, and $\{T^{(i)}\}$ satisfying \cref{equ:periodicKS} and the representation above.
\ENSURE
$\{Q^{(i)}\}_{i=1}^p$, $\{U_j^{(i)}\}_{i=1,j=1}^{p,d}$, and projected matrices $\{T^{(i)}\}_{i=1}^p$.
\vspace{0.25em}
  \FOR{$i=1,2,\ldots,p$}
    \STATE \textbf{(1) Generate new last-block vector $v\in\mathbb{R}^n$:}
    \STATE $v \leftarrow - (A_d^{(i)})^{-1}\sum_{\tau=1}^{d} A_{\tau-1}^{(i)}\, Q^{(i-d+\tau-1)}\, U_{\tau}^{(i-1)}(:,\text{end})$.
    \STATE \textbf{(2) First Level Orthonormalization:}
    \STATE $s \leftarrow (Q^{(i)})^{\mathrm{H}} v$;\quad $r_\perp \leftarrow v - Q^{(i)} s$;\quad $a \leftarrow \|r_\perp\|_2$;\quad $q\leftarrow r_\perp/a$;
    \STATE $Q^{(i)} \leftarrow [\,Q^{(i)},\, q\,]$;\quad $u_{v} \leftarrow [\,s;\,a\,]$.
    \STATE \textbf{(3) Second level Orthonormalization:}
    \STATE $u_{\tau} \leftarrow U_{\tau+1}^{(i-1)}(:,\text{end})$ for $\tau=1,\ldots,d-1$;\quad $u_d\leftarrow u_v$.
    \STATE $U_\tau^{(i)} \leftarrow \bigl[\,U_\tau^{(i)};\ 0\,\bigr]$ for $\tau=1,\dots,d,$ to conform with the augmented $Q^{(i-d+\tau)}$.
    \STATE $t \leftarrow \sum_{\tau=1}^{d} (U_\tau^{(i)})^{\mathrm{H}} u_{\tau}$.
    \STATE $u_\tau \leftarrow u_\tau - U_\tau^{(i)} t$ for $\tau=1,\ldots,d$; $\rho \leftarrow \sqrt{\sum_{\tau=1}^{d} \|u_\tau\|_2^2}$.
    \STATE $U_\tau^{(i)} \leftarrow [\,U_\tau^{(i)},\, u_\tau/\rho\,]$ for $\tau=1,\ldots,d$.
    \STATE \textbf{(4) Update projected matrix $T^{(i)}$:}
    \STATE $T^{(i)} \leftarrow \begin{bmatrix}
      T^{(i)} & t\\ 0 & \rho
    \end{bmatrix}.$
    \ENDFOR
\end{algorithmic}
\end{algorithm}

\subsection{Complexity Analysis}
We analyze the arithmetic and memory costs for large-scale LTP systems. Three workflows are compared: one-step collocation followed by \texttt{pKS}, multistep discretization followed by \texttt{pKS}, and multistep discretization followed by \texttt{pTOAR}. Only the dominant \(k\) Floquet multipliers are sought, with \(k\ll n\). 

Assume that each sampled coefficient matrix \(G(t_i)\) in \cref{equ:LTPOriginal} has \(\mathcal{O}(n)\) nonzeros and hence requires \(\mathcal{O}(n)\) storage. The symbol \(d\) has a workflow-dependent meaning: it denotes the number of steps, or history length, in the multistep workflows, and the number of collocation points per mesh interval in the collocation workflow. For collocation methods, the resulting one-step matrices are treated as dense after condensation. For the multistep discretizations, the matrices \(A_j^{(i)}\) inherit the sparsity of \(G(t_i)\), so each matrix--vector product costs \(\mathcal{O}(n)\). Let \(\mathcal C_{\mathrm{sol}}(n)\) denote the cost of one sparse solve with \(A_d^{(i)}\); the same sparse solver is used in both multistep workflows.

The one-step collocation+\texttt{pKS} workflow requires a condensation step to form \cref{equ:OneStepEVP}, with arithmetic cost \(\mathcal{O}(pn^3d^3)\) \cite{Lust2001improvedFloMul}. The resulting dense coefficient matrices require \(\mathcal{O}(pn^2)\) memory. Forming \(pk\) Arnoldi vectors costs \(\mathcal{O}(pn^2k)\), orthogonalization costs \(\mathcal{O}(pnk^2)\), and Ritz value extraction costs \(\mathcal{O}(pk^3)\). The periodic Arnoldi bases and projected matrices require \(\mathcal{O}(pnk)\) and \(\mathcal{O}(pk^2)\) memory, respectively.

For the multistep+\texttt{pKS} workflow, no condensation is needed. Forming each new Arnoldi vector requires \(d\) sparse matrix--vector products and one sparse solve. Over \(pk\) vectors, this costs \(\mathcal{O}\bigl(pk(nd+\mathcal C_{\mathrm{sol}})\bigr)\). Full-space orthogonalization adds \(\mathcal{O}(pndk^2)\), and Ritz value extraction adds \(\mathcal{O}(pk^3)\). The sparse coefficient matrices, full periodic Arnoldi bases, and projected matrices require \(\mathcal{O}(pn)\), \(\mathcal{O}(pndk)\), and \(\mathcal{O}(pk^2)\) memory, respectively.

For the multistep+\texttt{pTOAR} workflow, new-vector generation has the same cost,
\(\mathcal{O}\bigl(pk(nd+\mathcal C_{\mathrm{sol}})\bigr)\), as in multistep+\texttt{pKS}. Reconstructing the \(d\) history blocks of the Arnoldi vectors costs an additional \(\mathcal{O}(pndk^2)\). The first- and second-level orthogonalizations add \(\mathcal{O}(pnk^2)\) and \(\mathcal{O}(pdk^3)\), respectively, while Ritz value extraction costs \(\mathcal{O}(pk^3)\). The compressed bases and coefficient matrices require \(\mathcal{O}(pnk+pdk^2)\) memory, in addition to \(\mathcal{O}(pn)\) for the sparse coefficient matrices and \(\mathcal{O}(pk^2)\) for the projected matrices.

The componentwise estimates are summarized in \cref{tab:cost_comparison}. They are given for one complete, unrestarted periodic Arnoldi expansion generating \(pk\) Arnoldi vectors over one period; restart and any subsequent outer iterations are not included. The table reports the leading orders in \(n\), and the simplified total-memory bound for \texttt{pTOAR} uses \(dk\ll n\).
\begin{table}[htbp]
\centering
\small
\begin{tabular}{lccc}
\toprule
\textbf{Workflow} & \textbf{Collocation} & \textbf{Multistep}& \textbf{Multistep}\\
     & \textbf{+pKS} & \textbf{+pKS} &\textbf{+pTOAR}\\
\midrule
\multicolumn{4}{l}{\textit{Memory Cost}} \\
Coefficient Matrices  & $pn^{2}$ & $\mathcal{O}(pn)$ & $\mathcal{O}(pn)$  \\
Periodic Arnoldi Basis & $pnk$  & $pndk$ &$pnk+pdk^{2}$ \\
\textbf{Total Memory}           & $\mathcal{O}(pn^{2})$  & $\mathcal{O}(pndk)$ & $\mathcal{O}(pnk)$ \\
\midrule
\multicolumn{4}{l}{\textit{Arithmetic Cost}} \\
Condensation & $\mathcal{O}(pd^{3}n^3)$ & -- & -- \\
Arnoldi Vector Generation & $\mathcal{O}(pn^2k)$ & $\mathcal{O}(pk(nd+\mathcal C_{\mathrm{sol}}))$ & $\mathcal{O}(pk(nd+\mathcal C_{\mathrm{sol}}))$ \\
Orthogonalization & $\mathcal{O}(pnk^{2})$ & $\mathcal{O}(pndk^{2})$ & $\mathcal{O}(pnk^2+pdk^3)$ \\
Ritz Value Extraction & $\mathcal{O}(pk^3)$ & $\mathcal{O}(pk^3)$ & $\mathcal{O}(pk^3)$ \\
\textbf{Total Arithmetic Cost} & $\mathcal{O}(pn^{3}d^{3})$ & $\mathcal{O}(pk(ndk+\mathcal C_{\mathrm{sol}}))$ & $\mathcal{O}(pk(ndk+\mathcal C_{\mathrm{sol}}))$ \\
\bottomrule
\end{tabular}
\caption{Comparison of memory and arithmetic costs. The input matrices are assumed sparse with $\mathcal{O}(n)$ storage; $d$ denotes the number of steps of the multistep methods or the number of collocation points, and $\mathcal{C}_{\mathrm{sol}}=\mathcal{C}_{\mathrm{sol}}(n)$ denotes the cost of one sparse linear system solve.}
\label{tab:cost_comparison}
\end{table}

The collocation-based workflow is generally more expensive because condensation destroys sparsity. For an implicit \(d\)-step method, generating each new vector requires one sparse linear-system solve, independently of \(d\), together with \(d\) sparse matrix--vector products. Within a multistep family, increasing \(d\) typically permits a higher consistency order and hence higher accuracy. Since practical methods usually use \(d<5\), we have \(nd=\mathcal{O}(n)\); therefore, when \(\mathcal C_{\mathrm{sol}}(n)\gg n\), increasing \(d\) leaves the dominant solve cost unchanged and adds only lower-order matrix--vector work. This accuracy--cost advantage is shared by multistep+\texttt{pKS} and \texttt{pTOAR}, which have the same leading arithmetic complexity. \texttt{pTOAR} provides the additional advantage of reducing the Arnoldi-basis memory from \(\mathcal{O}(pndk)\) to \(\mathcal{O}(pnk+pdk^2)\). When \(dk\ll n\), the \(d\)-dependent term is lower order, so increasing the history length causes little growth in the dominant memory cost.

\section{Numerical Results}
\label{sec:results}
Three experiments assess the predicted convergence and parasitic decay, compare pTOAR with collocation-based packages (MATCONT, AUTO-07p), and demonstrate its performance on a large sparse RF system with discrete-time samples. The storage reduction of pTOAR follows from the compact representation and the complexity analysis in \cref{sec:Algorithm}; the experiments below assess convergence and runtime scaling. All computations use MATLAB R2024b on an Intel Core i5-1135G7 (2.4 GHz, 16 GB RAM).

To compare methods with different numbers of steps, we extract the last
\(n\)-dimensional block of each stacked invariant subspace, which corresponds to the common phase \(t_0\). Let \(e_{d}\) denote the $d$-th column of the identity matrix \(I_d\). For a \(d\)-step computation and a \(d_{\mathrm{ref}}\)-step reference, we define the relative error for Floquet multipliers and subspaces as
\[
e_{\mathrm{mul}}
=
\frac{|\widehat{\lambda}-\lambda_{\mathrm{ref}}|}
     {|\lambda_{\mathrm{ref}}|},
\qquad
e_{\mathrm{sub}}^{(m)}
=\left\lVert \sin\Theta\!\left(\bigl(e_d^{\mathrm{T}}\otimes I_n\bigr)\widehat{\mathcal{U}}_m,\bigl(e_{d_{\mathrm{ref}}}^{\mathrm{T}}\otimes I_n\bigr)\mathcal{U}_{m,\mathrm{ref}}
\right) \right\rVert.
\]
Here, \(\widehat{\mathcal{U}}_m\) is the computed stacked invariant subspace, and \(\mathcal{U}_{m,\mathrm{ref}}\) is its reference counterpart. When exact Floquet data are available, we take \(\lambda_{\mathrm{ref}}=\lambda\), \(d_{\mathrm{ref}}=d\), and \(\mathcal{U}_{m,\mathrm{ref}}=\mathcal{U}_m\), with \(\mathcal{U}_m\) defined in \cref{def:ExactStackedSelectedModes}. Otherwise, \(\lambda_{\mathrm{ref}}\) and \(\mathcal{U}_{m,\mathrm{ref}}\) are computed using BDF4 on the finest available grid, so that \(d_{\mathrm{ref}}=4\). For \(m=1\), \(e_{\mathrm{sub}}^{(1)}\) is the sine of the angle between the extracted computed and reference eigendirections.

\subsection{Verification of the Convergence Results}

The first experiment uses a two-dimensional LTP system \cref{equ:2DToyModelLTP} with analytically known Floquet eigenpairs.
\begin{equation}
    \begin{aligned}
        \begin{bmatrix}
            \dot{x}_{1}\\
            \dot{x}_{2}
        \end{bmatrix}=
        \begin{bmatrix}
            \alpha \cos^{2}t +(\beta-2)\sin^{2}t &1+(\beta-\alpha-2)\sin t\cos t\\
            -1+(\beta-\alpha-2)\sin t\cos t& \alpha \sin^{2}t +(\beta-2) \cos^{2}t 
            \end{bmatrix}\begin{bmatrix}
            x_1\\x_2
        \end{bmatrix}.
    \end{aligned}
    \label{equ:2DToyModelLTP}
\end{equation}
Its Floquet multipliers are $e^{2\pi \alpha}$ and $e^{2\pi(\beta-2)}$, and their corresponding eigenvectors are columns of 
\begin{equation*}
    \begin{aligned}
        U(t) = \begin{bmatrix}
             \cos t& \sin t\\-\sin t& \cos t
        \end{bmatrix}.
    \end{aligned}
\end{equation*}
In this experiment, we fix $\alpha=\beta=0.1$.
One period $[0,2\pi]$ is discretized into $2^\ell+1$ grid points using a non-uniform mesh. To reflect a general variable stepsize setting while satisfying Assumption~(A2), two local stepsize enlargements are introduced near $t=\pi/2$ and two reductions near $t=3\pi/2$ within each period.

BDF2 and BDF3 are applied to form the pPEP \cref{equ:pPEPDefinition}. All Floquet multipliers, eigenvectors, and parasitic periodic eigenvalues are computed by the periodic QR algorithm \cite{bojanczyk1992Periodic}. As a one-step comparator, we also solve the BE-discretized periodic eigenvalue problem; unlike the multistep pPEPs, it has no parasitic periodic eigenvalues.

\begin{figure}[!htbp]  
   \centering
   \begin{subfigure}[b]{0.48\linewidth}
        \centering
        \includegraphics[width=0.8\linewidth]{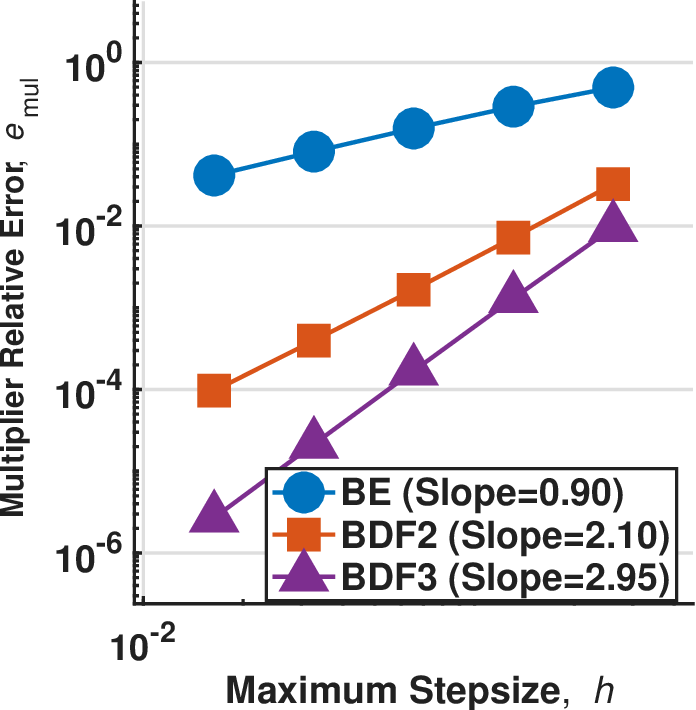}
        \caption{Floquet Multiplier Error vs. Maximum Stepsize}
        \label{fig:EigValErr}
    \end{subfigure}
    \begin{subfigure}[b]{0.48\linewidth}
        \centering
        \includegraphics[width=0.8\linewidth]{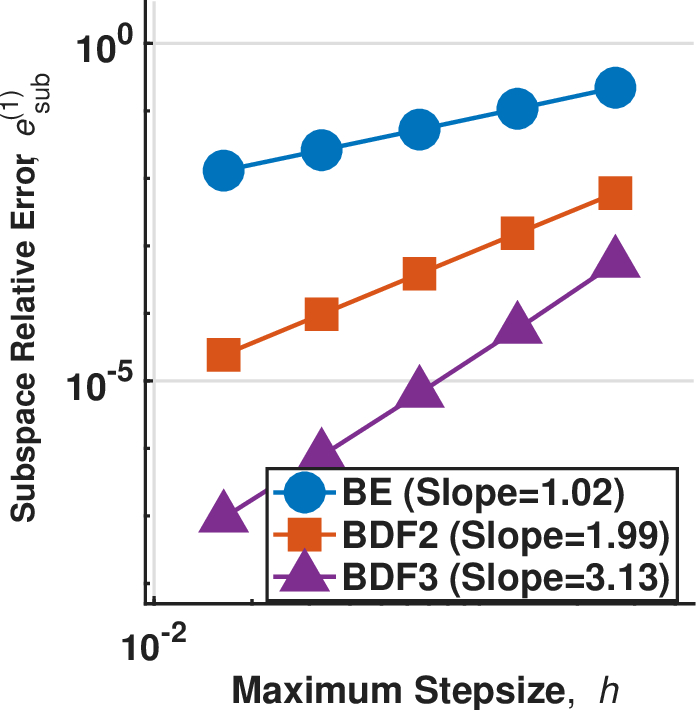}
        \caption{Floquet Eigenvector Error vs. Maximum Stepsize}
        \label{fig:EigVecErr}
    \end{subfigure}
    \caption{Convergence Behavior of Dominant Floquet Eigenpair.}
    \label{fig:SL} 
\end{figure}

Relative errors for the dominant Floquet multiplier and its eigenvector are shown in \cref{fig:EigValErr} and \cref{fig:EigVecErr}. BE, BDF2, and BDF3 exhibit first-, second-, and third-order convergence, respectively. The BDF2 and BDF3 results are consistent with \cref{thm:MainConvergence}; BElies outside its $d\ge2$ setting. Results for the other eigenpair are similar.

\begin{figure}[!htbp]
    \centering
        \includegraphics[width=0.6\linewidth]{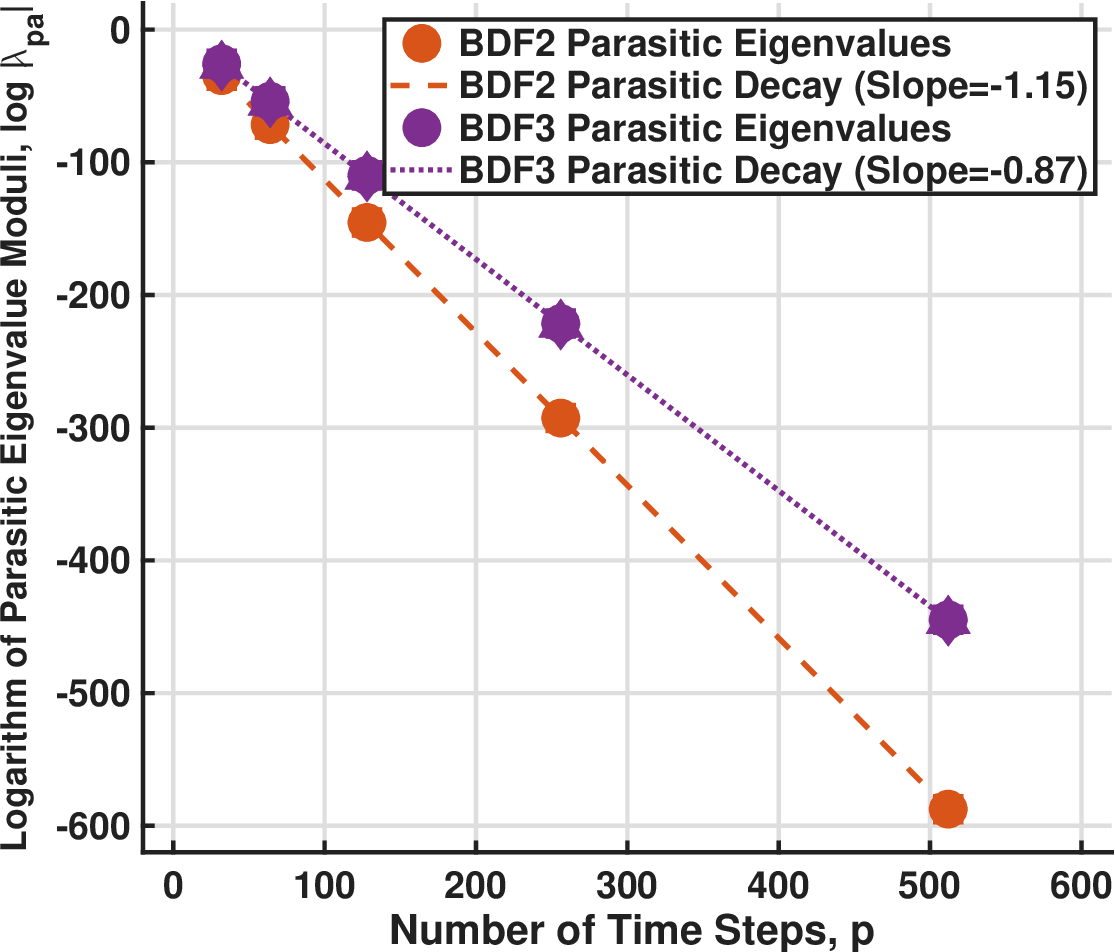}
    \caption{Convergence Behavior of Parasitic Periodic Eigenvalues.}
    \label{fig:SLParasitic} 
\end{figure}

\Cref{fig:SLParasitic} shows the logarithmic magnitudes of the parasitic periodic eigenvalues of BDF2 and BDF3 as the total number of time steps increases. For each discretization level, the parasitic eigenvalues have nearly identical magnitudes; hence their logarithmic values visually coincide in \cref{fig:SLParasitic}. The trend lines are obtained by fitting the averaged logarithmic magnitudes at each step number. For reference, the theoretical parasitic roots are \(\nu_{2}=\tfrac{1}{3}\) for BDF2, and \(\nu_{2,3}=\frac{7\pm \jmath\sqrt{39}}{22}\) for BDF3, leading to expected slopes \(\log|\nu_{2}|\approx -1.10\) and \(\log|\nu_{2,3}|\approx -0.85\), respectively. The fitted slopes (about \(-1.15\) for BDF2 and \(-0.87\) for BDF3) are consistent with the \(|\nu_\tau|^p\) decay predicted by \cref{thm:MainConvergence}.
 
\subsection{Application 1: Coupled Stuart-Landau Oscillators}

We next consider a parameter-dependent coupled Stuart--Landau system of modest dimension, adapted from \cite{OK1991Fairgrieve, Lust2001improvedFloMul}. For this application, the key numerical requirement is to reliably classify each Floquet multiplier as lying inside, outside, or on the unit circle. This classification determines the stability of the periodic orbit and the dimension of its unstable manifold \cite{Lust2001improvedFloMul}; high accuracy of each individual multiplier is not required.

\begin{equation*}
    \left\{
    \begin{aligned}
        \dot{x}_{1} &= x_{1}+\beta y_{1} - x_{1}(x_{1}^{2}+y_{1}^{2}) + \delta (x_{2}-x_{1}+y_{2}-y_{1}),\\
        \dot{y}_{1} &= -\beta x_{1}+ y_{1} - y_{1}(x_{1}^{2}+y_{1}^{2}) + \delta (x_{2}-x_{1}+y_{2}-y_{1}),\\
        \dot{x}_{2} &= x_{2}+\beta y_{2} - x_{2}(x_{2}^{2}+y_{2}^{2}) - \delta (x_{2}-x_{1}+y_{2}-y_{1}),\\
        \dot{y}_{2} &= -\beta x_{2}+ y_{2} - y_{2}(x_{2}^{2}+y_{2}^{2}) - \delta (x_{2}-x_{1}+y_{2}-y_{1}).\\
    \end{aligned}\right.
\end{equation*} 

We fix $\beta=0.5$ and use $\delta$ as a continuation parameter, approaching $\delta=0.25$ from below. For each $\delta$, we consider a limit cycle satisfying \(x_1 = -x_2, y_1 = -y_2 \). As $\delta \to 0.25$, this limit cycle approaches a heteroclinic loop, and the period goes to infinity \cite{Lust2001improvedFloMul}. For each $(\delta,\beta)$, the periodic orbit is computed using continuation packages (\texttt{MATCONT} and \texttt{AUTO-07p}), which also provide the Jacobian matrices along the trajectory, i.e., the samples of the variational LTP system's coefficient matrices.

We compare three discretization-solver workflows for computing the Floquet multipliers of this variational system. \texttt{AUTO-07p} employs collocation together with the O.K. Floquet algorithm \cite{OK1991Fairgrieve}. A modified version of \texttt{MATCONT} retains the collocation discretization but replaces the built-in MATLAB \texttt{eig} routine with a periodic QR algorithm \cite{bojanczyk1992Periodic}. Our proposed approach discretizes the variational system using BDF4 and computes multipliers via \texttt{pTOAR}. For a fair comparison, each workflow uses a total of 1024 function evaluations per period. For the direct methods \texttt{AUTO-07p} and the modified \texttt{MATCONT}, we set \texttt{NTST}=512 and \texttt{NCOL}=2; both methods compute all four Floquet multipliers. For the iterative BDF4--\texttt{pTOAR} workflow, we compute the four principal Floquet multipliers in a single run; the parasitic periodic eigenvalues are not computed.

\begin{figure}[!htbp]
    \centering
    \begin{subfigure}[b]{0.48\linewidth}
        \centering
        \includegraphics[width=\linewidth]{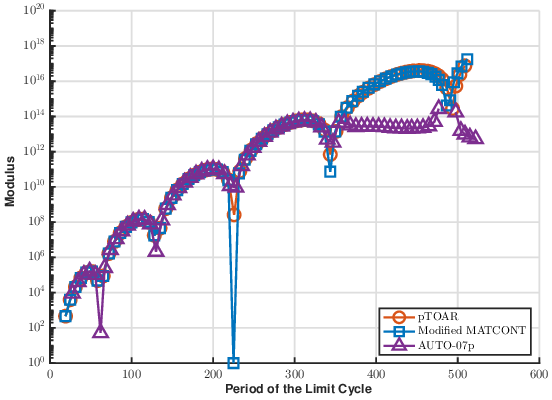}
        \caption{First Floquet Multiplier}
    \end{subfigure}
    \begin{subfigure}[b]{0.48\linewidth}
        \centering
        \includegraphics[width=\linewidth]{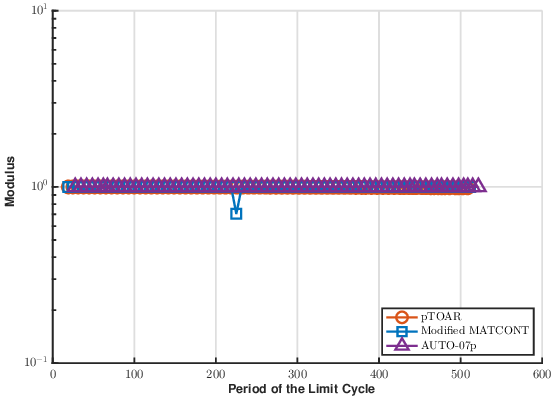}
        \caption{Second Floquet Multiplier}
    \end{subfigure}

    \begin{subfigure}[b]{0.48\linewidth}
        \centering
        \includegraphics[width=\linewidth]{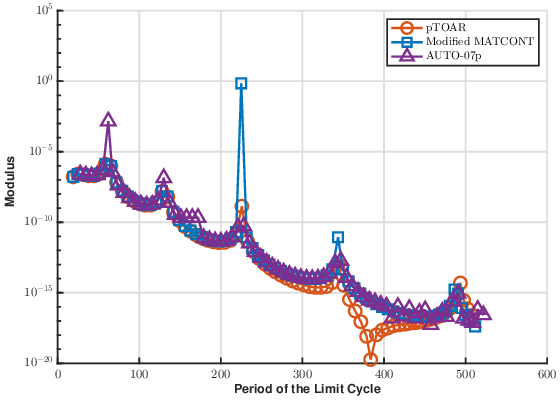}
        \caption{Third Floquet Multiplier}
    \end{subfigure}
    \begin{subfigure}[b]{0.48\linewidth}
        \centering 
        \includegraphics[width=\linewidth]{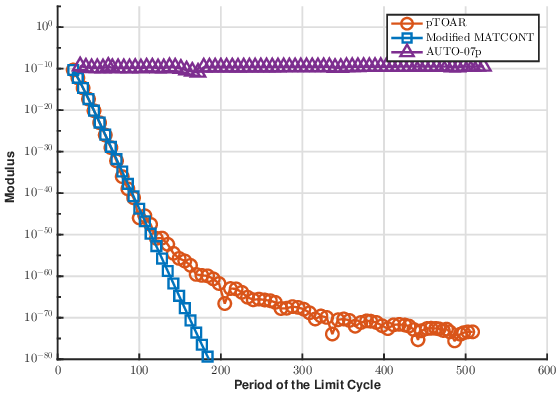}
        \caption{Fourth Floquet Multiplier}
    \end{subfigure}
    \caption{Evolution of Floquet Multipliers vs. Period of the Parameterized Limit Cycle.}
    \label{fig:Brus_Flomul1234} 
\end{figure}

The moduli of the computed Floquet multipliers are plotted on a logarithmic scale versus the period in \cref{fig:Brus_Flomul1234}; sharp transitions near suspected bifurcation points are not fully resolved at the computed parameter values.

All three methods yield nearly identical curves for the three largest multipliers, except that \texttt{AUTO-07p} fails to report the largest one once it exceeds $10^{16}$. For the smallest multiplier, both the modified \texttt{MATCONT} and our method show a monotonic decay, whereas \texttt{AUTO-07p} reports an approximately constant value around $10^{-10}$ over a wide range of periods.

\begin{figure}[!htbp]
    \centering
        \includegraphics[width=0.7\linewidth]{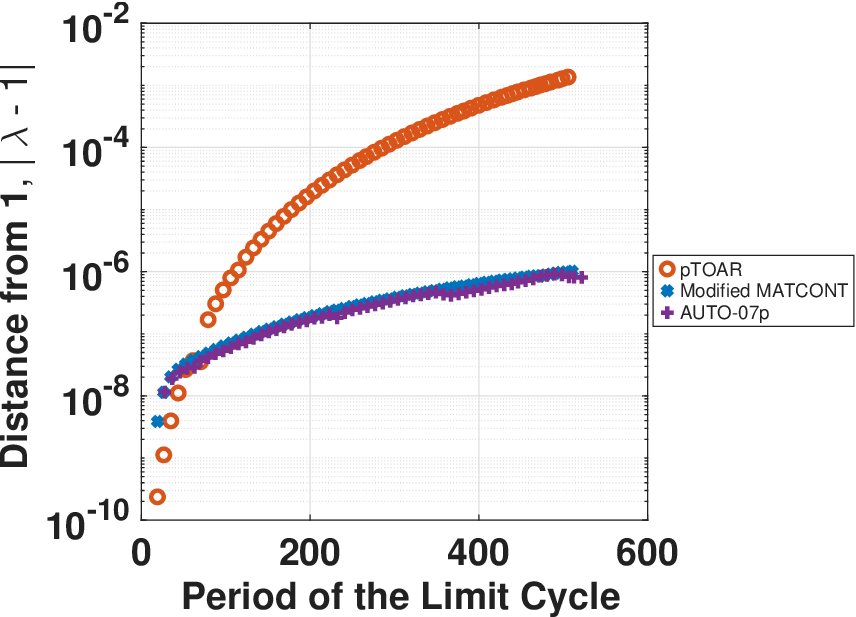}
    \caption{Distance of the Oscillatory Floquet Multiplier from 1 vs. Period of the Parameterized Limit Cycle.}
    \label{fig:Err_Brus_Flomul} 
\end{figure} 

\Cref{fig:Err_Brus_Flomul} shows $|\lambda-1|$ for the oscillatory Floquet multiplier, whose exact value is $1$. \texttt{AUTO-07p} and the modified \texttt{MATCONT} maintain high accuracy with fourth-order collocation. Under the fixed budget of 1024 function evaluations per period, BDF4 is less accurate for longer periods. Nevertheless, it keeps the oscillatory multiplier close to $1$ and provides the correct stability classification throughout the tested period range.

\subsection{Application 2: RF Circuits}

In this experiment, we consider a large-scale, sparse LTP system arising from radio-frequency (RF) circuit simulation, for which the system matrices are available only at discrete time samples. The target quantity is the perturbation projection vector (PPV), which is used in phase-noise analysis \cite{demir2000floquet}. The PPV is a nontrivial $T$-periodic solution of the adjoint circuit DAE
\begin{equation}
-\begin{bmatrix}
    C_{11}(t) & 0\\
    0 & 0
\end{bmatrix}^{\mathrm{T}}
\frac{\mathrm{d}}{\mathrm{d}t}
\begin{bmatrix}
    y_{1}\\
    y_{2}
\end{bmatrix}
+
\begin{bmatrix}
    G_{11}(t) & G_{12}(t)\\
    G_{21}(t) & G_{22}(t)
\end{bmatrix}^{\mathrm{T}}
\begin{bmatrix}
    y_{1}\\
    y_{2}
\end{bmatrix}
=0.
\label{equ:AdjointCircuitDAE}
\end{equation}
The coefficient matrices are $T$-periodic, and $C_{11}(t)$ and $G_{22}(t)$ are assumed to be nonsingular for all $t$. Thus, \cref{equ:AdjointCircuitDAE} is a semi-explicit index-1 DAE \cite{marz1995DAETheory}. Solving its algebraic block for $y_2$ and substituting the result into the differential block gives
\begin{equation}
\left\{
\begin{aligned}
&-\frac{\mathrm{d}}{\mathrm{d}t}y_{1}(t)
+C_{11}^{-\mathrm{T}}(t)
\left[
G_{11}^{\mathrm{T}}(t)
-G_{21}^{\mathrm{T}}(t)
G_{22}^{-\mathrm{T}}(t)
G_{12}^{\mathrm{T}}(t)
\right]y_{1}(t)=0,\\
&y_{2}(t)
=-G_{22}^{-\mathrm{T}}(t)G_{12}^{\mathrm{T}}(t)y_{1}(t).
\end{aligned}
\right.
\label{equ:DAEdecoupled}
\end{equation}

For a practical RF oscillator, the PPV is associated with the trivial Floquet multiplier at $1$, whereas all other physical Floquet multipliers lie inside the unit circle. We therefore compute the PPV as the periodic eigenvector associated with the computed multiplier closest to $1$.

We consider two RF circuit cases, referred to below as Case 1 and Case 2. The first case has a single isolated multiplier closest to $1$, for which eigenvector convergence is expected. The second case exhibits a cluster of multipliers near $1$, so that individual eigenvectors are ill-conditioned, and we assess convergence in terms of the associated invariant subspace. Basic case information is reported in \cref{tab:BasicInfoPPV}.
\begin{table}[h]
    \centering
    \begin{tabular}{cccc}
        \toprule[1.5pt]
        Model & Dimension & Differential dimension & Number of time steps\\
        \midrule[1pt]
        Case 1 &18 & 16 & 448\\
        Case 2 &913 & 388 &654\\
        \bottomrule[1.5pt]
    \end{tabular}
    \caption{Basic Information of RF Test Cases.}
    \label{tab:BasicInfoPPV}
\end{table}

The dimension in \cref{tab:BasicInfoPPV} refers to the original DAE \cref{equ:AdjointCircuitDAE}, while the differential dimension is the size of \(C_{11}(t)\) and hence the dimension of \cref{equ:DAEdecoupled}.

We use BDF4 on the full available grids of 448 and 654 samples as the numerical reference and apply BE, BDF2, and BDF3 with \texttt{pTOAR} to compute the Floquet eigenpairs. Coarser grids are formed by subsampling these prescribed grids; no interpolation is used. Because collocation would require values at additional internal points or an interpolation model, it is not used in this comparison.

\begin{figure}[!htbp]
    \centering
    \begin{subfigure}[b]{0.3\linewidth}
        \centering
        \includegraphics[width=\linewidth]{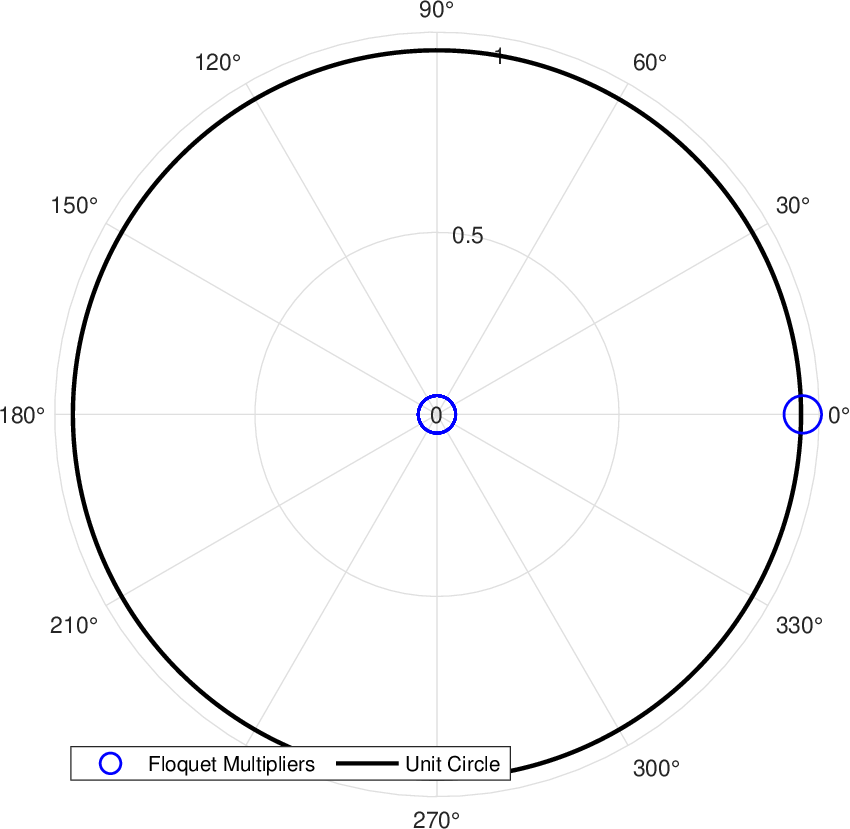}
        \caption{Floquet Multipliers Relative to the Unit Circle}
        \label{fig:RFCase1Polar}
    \end{subfigure}
    \hfill 
    \begin{subfigure}[b]{0.3\linewidth}
        \centering
        \includegraphics[width=\linewidth]{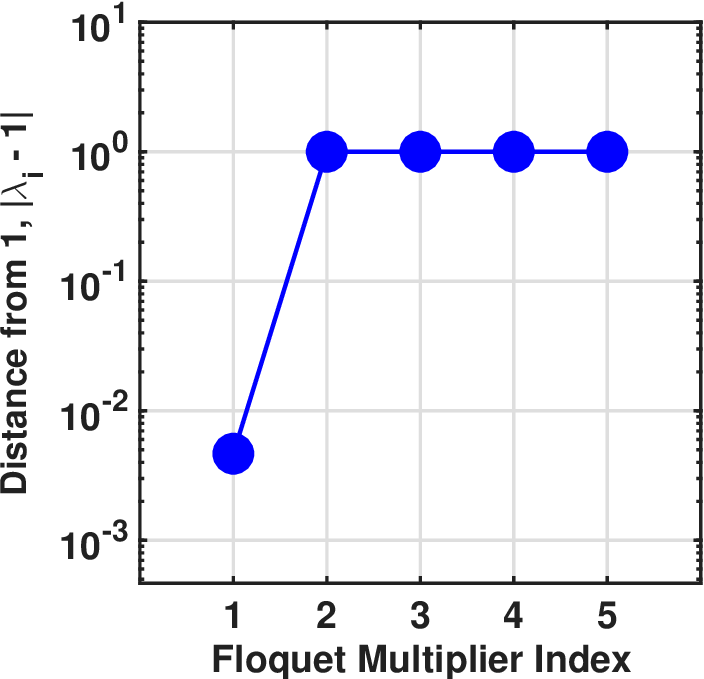}
        \caption{Distances between the Dominant Floquet Multipliers and 1}
        \label{fig:RFCase1Separation}
    \end{subfigure}
    \hfill 
    \begin{subfigure}[b]{0.3\linewidth}
        \centering
        \includegraphics[width=\linewidth]{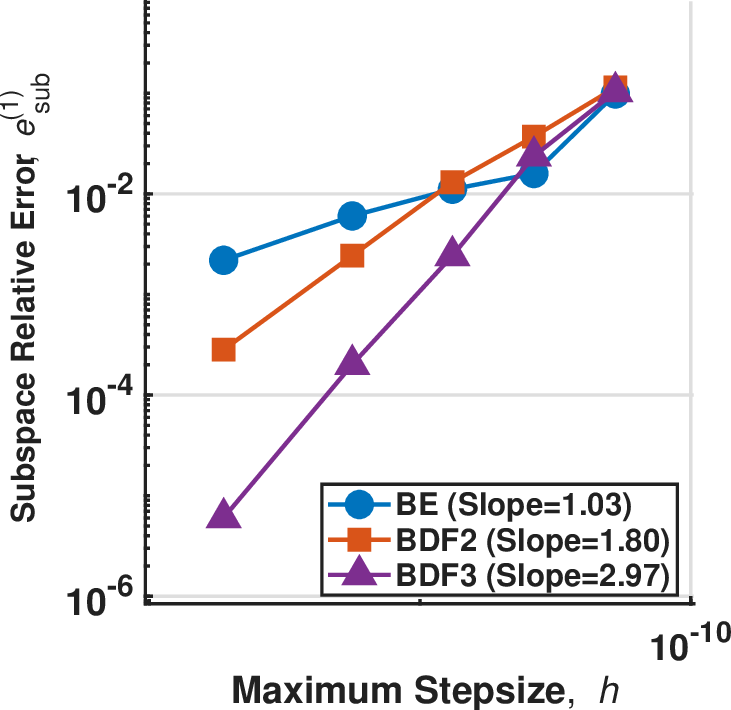} 
        \caption{Dominant Eigenvector Error vs. Maximum Stepsize}
        \label{fig:RFCase1Eigenvector}
    \end{subfigure}
    \caption{Floquet Multiplier Distribution and Dominant Eigenvector Convergence for Case 1.}
    \label{fig:RFCase1}
\end{figure}

As shown in \cref{fig:RFCase1Polar,fig:RFCase1Separation}, all Floquet multipliers except one are clustered near the origin; a single dominant multiplier lies close to $1$. By \cref{thm:MainConvergence}, this separation favors accurate recovery of its one-dimensional invariant subspace for BDF2 and BDF3; BE again serves as a one-step comparator. The results in \cref{fig:RFCase1Eigenvector} exhibit the expected first-, second-, and third-order convergence rates for BE, BDF2, and BDF3, respectively.

\begin{figure}[!htbp]
    \centering
    \begin{subfigure}[b]{0.48\linewidth}
    \centering
    \includegraphics[width=.8\linewidth]{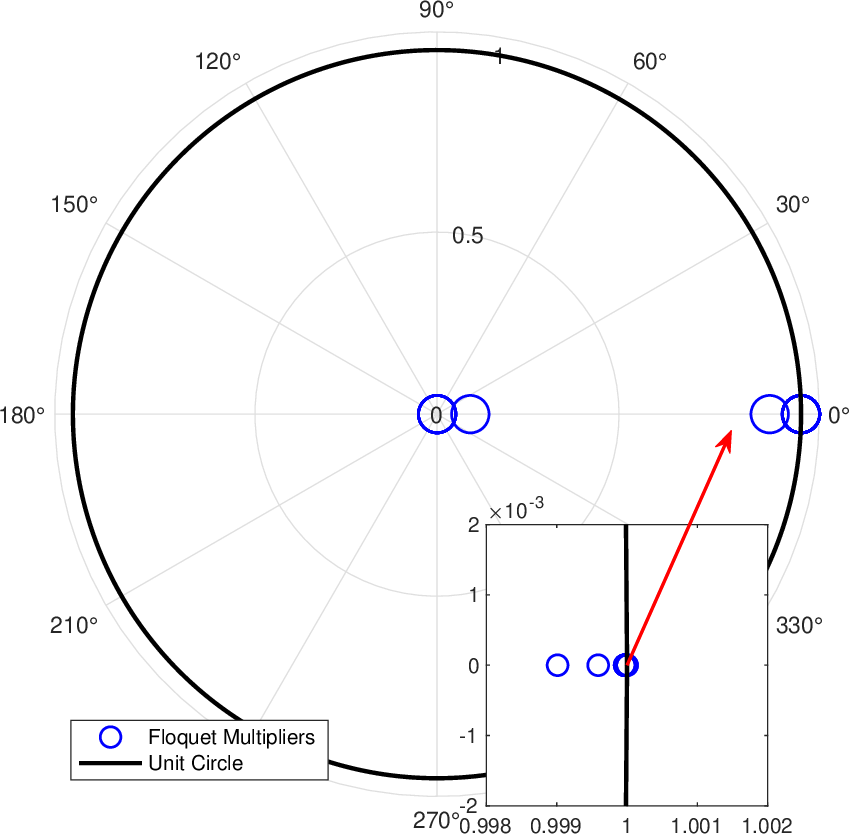}
    \caption{Floquet Multipliers Relative to the Unit Circle}
    \label{fig:RFCase2Polar}
    \end{subfigure}
    \hfill
    \begin{subfigure}[b]{0.48\linewidth}
    \centering
    \includegraphics[width=0.8\linewidth]{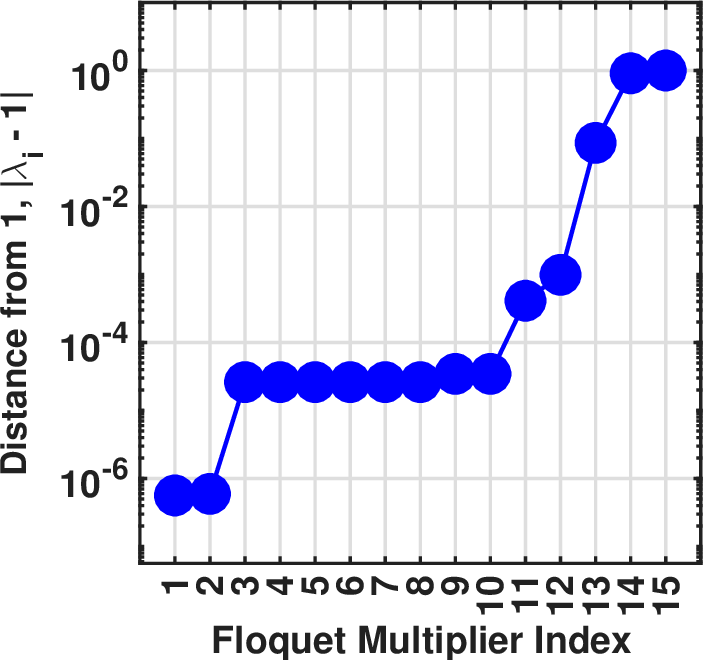}
    \caption{Distances between the Dominant Floquet Multipliers and 1}
    \label{fig:RFCase2Separation}
    \end{subfigure}
    \caption{Floquet Multiplier Distribution for Case 2.}
    \label{fig:RFCase2Spectrum}
\end{figure}

In contrast to Case 1, Case 2 has 13 multipliers tightly clustered near $1$ but separated from the remainder of the spectrum, as shown in \cref{fig:RFCase2Polar,fig:RFCase2Separation}. \Cref{thm:MainConvergence} predicts convergence of the 13-dimensional invariant subspace, but not of individual eigenvectors.

\begin{figure}[!htbp]
    \centering
    \begin{subfigure}[b]{0.48\linewidth}
        \centering
        \includegraphics[width=.8\linewidth]{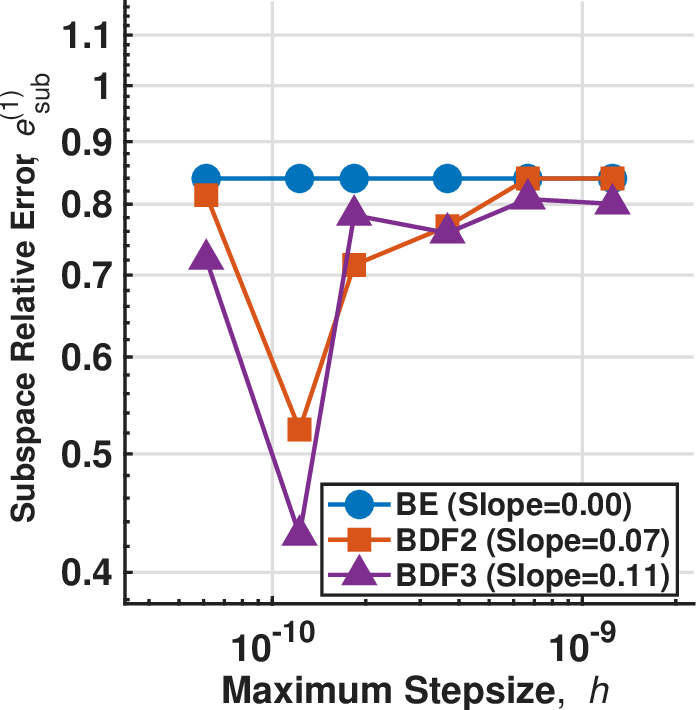}
        \caption{Dominant Eigenvector Error vs. Maximum Stepsize}
        \label{fig:RFCase2Eigenvector}
    \end{subfigure}
    \hfill
    \begin{subfigure}[b]{0.48\linewidth}
        \centering
        \includegraphics[width=.8\linewidth]{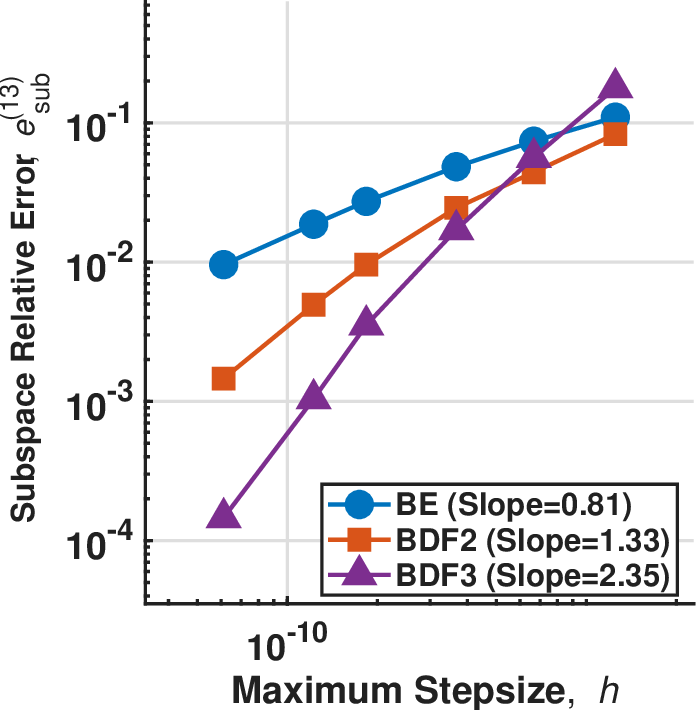}
        \caption{Dominant 13-Dimensional Subspace Error vs. Maximum Stepsize}
        \label{fig:RFCase2Subspace}
    \end{subfigure}
    \caption{Convergence Behavior of the Dominant Invariant Subspace in Case 2.}
    \label{fig:RFCase2Convergence}
\end{figure}

\Cref{fig:RFCase2Convergence} is consistent with this prediction over the tested grid range: the dominant Floquet eigenvector fails to converge, whereas the 13-dimensional invariant subspace exhibits the expected convergence behavior. The difficulty in recovering an individual PPV is therefore attribute to spectral ill-conditioning.

\begin{figure}[!htbp]
    \centering
    \begin{subfigure}[b]{0.48\linewidth}
        \centering 
        \includegraphics[width=.8\linewidth]{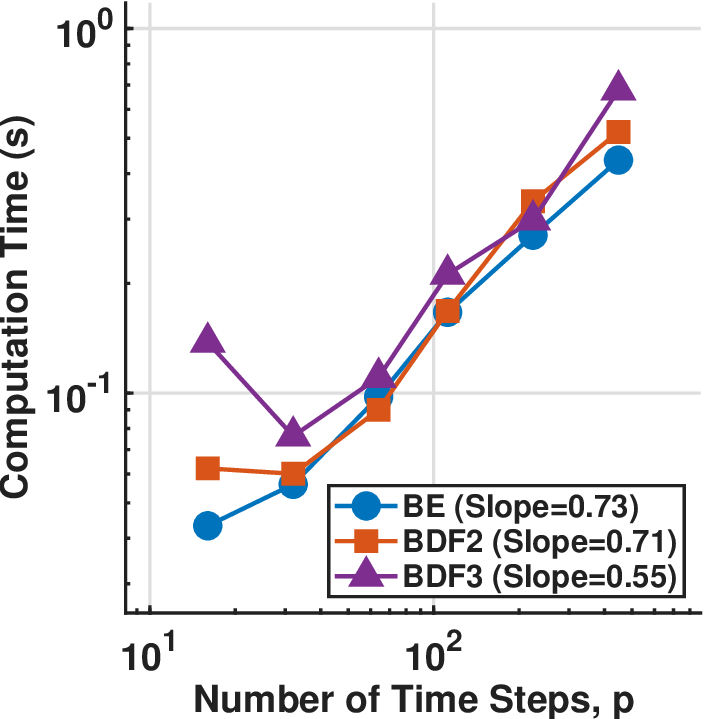}
        \caption{Case 1}
        \label{fig:RFCase1Timing}
    \end{subfigure}
    \hfill
    \begin{subfigure}[b]{0.48\linewidth}
        \centering
        \includegraphics[width=.8\linewidth]{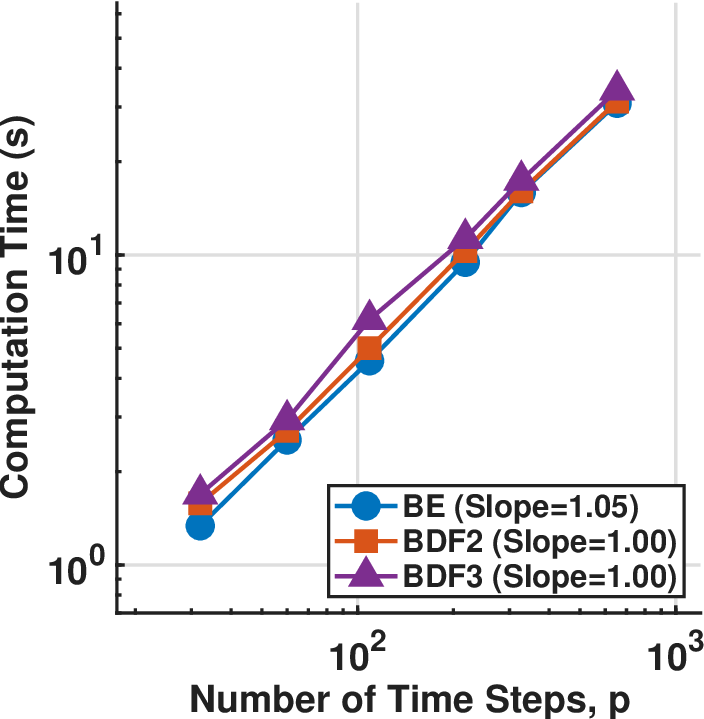}
        \caption{Case 2}
        \label{fig:RFCase2Timing} 
    \end{subfigure} 
    \caption{Computation Time vs. the Number of Time Steps.}
    \label{fig:ErrTiming}
\end{figure}

\cref{fig:ErrTiming} shows that the runtime grows linearly with the number of time steps and changes little with the choice of the multistep method. This agrees with the arithmetic complexity analysis in \cref{sec:Algorithm}. For Case 1, the tested problem sizes are too small for the asymptotic linear scaling to dominate, and the measured runtimes are affected by non-negligible per-run overhead.

\section{Conclusion}
\label{sec:conclu}
A multistep discretization of an LTP system gives rise to a pPEP for computing Floquet multipliers and invariant subspaces. As the stepsize tends to zero, the principal periodic eigenvalues and selected spectrally separated invariant subspaces converge at the consistency order, and the parasitic periodic eigenvalues decay geometrically. The resulting spectral separation motivates computing only the dominant spectrum. The proposed \texttt{pTOAR} algorithm exploits the companion structure to preserve the leading arithmetic order of \texttt{pKS} applied to the linearization while avoiding its $d$-fold basis-storage growth. Numerical experiments exhibit the predicted convergence and parasitic-decay behavior and demonstrate the method on dynamical-system and RF-circuit problems.

\section*{Acknowledgments}
This work was supported by the National Key R\&D Program of China under Grant Nos. 2020YFA0711900 and 2020YFA0711902. We thank our industrial partners for providing the RF simulation cases used in this study. We are grateful to the reviewers for their constructive comments and suggestions, and to Kimi for its assistance with language polishing.

\bibliographystyle{siamplain}
\bibliography{references}
\end{document}